# FORMES MODULAIRES $p$-ADIQUES

*par*

Antoine Chambert-Loir

**Résumé**. — Ce texte correspond à un exposé fait lors de la semaine «Formes modulaires et représentations galoisiennes» (Luminy, 1997). On expose la théorie des formes modulaires $p$-adiques, développée par Serre, Katz, Hida, Wiles, Coleman...

**Abstract**. — This is the text of a talk to the study week on *Modular forms and Galois representations* held in Luminy, 1997. We give a survey of $p$-adic modular forms, as developped by Serre, Katz, Hida, Wiles, Coleman and others...

## Table des matières

## Introduction

La motivation initiale de la théorie des formes modulaires $p$-adiques était de construire la fonction zêta $p$-adique d'un corps totalement réel quelconque. Rappelons simplement le cas de $\mathbf{Q}$. Considérons la fonction zêta de Riemann définie pour $\operatorname{Re}(s) > 1$ par $\zeta(s) = \sum_{n=1}^{\infty} n^{-s}$ ; elle admet un prolongement méromorphe à $\mathbf{C}$, avec un unique pôle (simple) en $s = 1$. Si $n$ est un entier $\geqslant 1$, on a la formule $\zeta(1-n) = -B_n/n$, le nombre de Bernoulli $B_n$ étant défini par le développement

$$\frac{z}{e^z - 1} = \sum_{n=0}^{\infty} B_n \frac{z^n}{n!}.$$

(Ils sont nuls pour $n$ impair $\geqslant 3$.) En particulier, les valeurs de la fonction zêta aux entiers négatifs sont des nombres rationnels. Kubota et Leopoldt ont découvert que ces valeurs sont $p$-adiquement « continues » et ont construit des fonctions analytiques $\zeta_{p,\alpha} : \mathbf{Z}_p \to \mathbf{Z}_p$ (resp. $\mathbf{Z}_p \setminus \{1\} \to \mathbf{Z}_p$ si $\alpha = 0$, avec un pôle simple en $s = 1$), indexées par $\alpha \in 2\mathbf{Z}/(p-1)\mathbf{Z}$, telles que pour tout entier pair $n \geqslant 2$,

$$\zeta_{p,\alpha}(1-n) = (1 - p^{n-1})\zeta(1-n), \qquad n \equiv \alpha \pmod{(p-1)}.$$

Pour $k$ pair, la série d'Eisenstein

$$G_k = -\frac{B_k}{2k} + \sum_{n=1}^{\infty} \sigma_{k-1}(n)q^n = -\frac{1}{2}\zeta(1-k) + \sum_{n=1}^{\infty} \sigma_{k-1}(n)q^n$$

est une forme modulaire (non parabolique) de poids $k$ pour le groupe $\mathrm{SL}(2,\mathbf{Z})$. (On note classiquement $\sigma_{k-1}(n) = \sum_{d|n} d^{k-1}$.) Alors, si $\sigma^*_{k-1}(n)$ est la somme des $d^{k-1}$ pour les diviseurs $d$ de $n$ qui sont premiers à $p$,

$$G_k^*(q) \stackrel{\mathrm{def}}{=} G_k(q) - p^{k-1}G_k(q^p) = -\frac{1}{2}(1-p^{k-1})\zeta(1-k) + \sum_{n=1}^{\infty} \sigma^*_{k-1}(n)q^n$$

est une forme modulaire de poids $k$ de niveau $p$. Si $(k_i)$ est une suite d'entiers pairs tendant vers $+\infty$ dans $\mathbf{R}$ et vers $k$ dans $\mathbf{Z}_p$, il est aisé de voir que $a_n(G_{k_i})$ converge vers $a_n(G_k^*)$ pour tout $n \geqslant 1$ et ce, *uniformément*. La théorie des formes modulaires $p$-adiques développée par J-P. Serre lui a permis d'en déduire (dans le cas beaucoup plus général d'un corps de nombres totalement réel) la continuité $p$-adique du terme constant, et donc l'interpolation $p$-adique des valeurs aux entiers négatifs de la fonction zêta de Riemann.

Une première façon de définir les formes modulaires $p$-adiques est ainsi de considérer des séries $f \in \mathbf{Z}_p[[q]]$ qui sont $p$-adiquement limite uniforme de développements de Fourier attachés à des formes modulaires usuelles. Ce point de vue sera l'objet du § 1.

On désire, en introduisant les formes modulaires $p$-adiques, expliquer les propriétés purement $p$-adiques des formes modulaires, et en particulier les congruences $p$-adiques auxquelles elles donnent lieu. Si la théorie des formes modulaires a un sens sur tout anneau, la considération des formes modulaires sur $\mathbf{Z}_p$ ne peut pas fournir une réponse adéquate à notre question, puisque celles-ci ne seront autres que les formes modulaires sur $\mathbf{Z}$ tensorisées par $\mathbf{Z}_p$.

On veut par exemple qu'une congruence du type $E_{p-1} \equiv 1 \pmod{p}$ se reflète en une égalité $E_{p-1} = 1 + pf$, où $f$ serait une forme modulaire ($p$-adique). Or une telle égalité est impossible avec une forme modulaire usuelle $f$ : en évaluant l'égalité précédente en une courbe elliptique supersingulière sur $\mathbf{F}_p$, on obtient $0 = 1$, la série d'Eisenstein $E_{p-1}$ étant congrue à l'invariant de Hasse modulo $p$. Élémentairement, on pourrait se contenter de la série

$$f = (E_{p-1} - 1)/p = -\frac{2(p-1)}{pB_{p-1}} \sum_{n=1}^{\infty} \sigma_{p-2}(n)q^n$$

(rappelons que $pB_{p-1}$ est une unité $p$-adique). Il est en tout cas nécessaire « d'éviter » les courbes elliptiques supersingulières, et c'est ce que fait explicitement N. Katz en développant une théorie « modulaire » des formes modulaires $p$-adiques, *cf.* le § 2.

Nous exposerons ensuite au § 3 l'algèbre de Hecke $p$-adique, construite par H. Hida. Si $\Lambda = \mathbf{Z}_p[[1 + p\mathbf{Z}_p]] \simeq \mathbf{Z}_p[[T]]$, c'est une $\Lambda$-algèbre. Nous définirons les formes modulaires *ordinaires*, ainsi que le facteur direct de l'algèbre de Hecke $p$-adique qui leur est attaché et conclurons ce paragraphe en énonçant les théorèmes de Hida sur la structure de cette algèbre de Hecke ordinaire.

Le paragraphe § 4 est consacré aux formes modulaires $\Lambda$-adiques. Convenons ici qu'une forme modulaire $\Lambda$-adique est une série

$$F = \sum_{n=0}^{\infty} a_n(T)q^n, \qquad q_n \in \Lambda = \mathbf{Z}_p[[T]]$$

telle que ses *spécialisations*

$$\sum_{n=0}^{\infty} a_n((1+p)^k - 1)q^n$$

sont pour (presque) tout $k \geqslant 2$ le $q$-développement d'une forme modulaire (usuelle) $f_k$ de poids $k$.

Nous expliquerons alors comment A. Wiles, redémontrant un théorème de Hida[1], attache à une forme modulaire $\Lambda$-adique une représentation galoisienne

$$\rho : \mathrm{Gal}(\overline{\mathbf{Q}}/\mathbf{Q}) \to \mathrm{GL}(2, \Lambda)$$

qui par la spécialisation $T \to (1+p)^k - 1$ redonne la représentation galoisienne attachée par Deligne à la forme modulaire $f_k$.

Nous terminerons cet exposé en décrivant au § 5 les formes modulaires $p$-adiques surconvergentes introduites par Katz, et revues du point de vue rigide-analytique par R. Coleman.

*Je remercie Jacques Tilouine pour ses conseils lors de la préparation de cet exposé.*

## Notations et conventions

Dans tout l'exposé, on fixe un nombre premier $p \geqslant 5$. Soit $\mathbf{C}_p$ la complétion d'une clôture algébrique de $\mathbf{Q}_p$. On fixe deux plongements $\overline{\mathbf{Q}} \hookrightarrow \mathbf{C}$ et $\overline{\mathbf{Q}} \hookrightarrow \mathbf{C}_p$. La norme $p$-adique d'un élément de $\mathbf{C}_p$ est normalisée par $|p|_p = 1/p$.

---

[1] en le généralisant, car le résultat de Wiles est aussi valable pour les formes modulaires ordinaires d'un corps totalement réel quelconque.

On appelle anneau $p$-adique tout anneau $A$ tel que $A = \varprojlim_n A/p^n A$.

Si $N$ est un entier $\geqslant 1$, on définit trois sous-groupes $\Gamma_0(N) \supset \Gamma_1(N) \subset \Gamma(N)$ de $\mathrm{SL}(2, \mathbf{Z})$ dont les éléments sont les matrices $\begin{pmatrix} a & b \\ c & d \end{pmatrix}$ telles que $c$ est multiple de $N$ (resp. de plus, $d \equiv 1 \pmod{N}$, resp. de plus, $b \equiv 0 \pmod{N}$).

## 1. Développements de Fourier

Dans ce paragraphe, les formes modulaires $p$-adiques que l'on introduit sont des limites $p$-adiques de développements de Fourier de formes modulaires usuelles. C'est ainsi que Serre [25] les a définies originellement, lorsque le niveau est 1. Néanmoins, l'extension de la théorie à un niveau quelconque, due à Katz, nécessite l'approche plus sophistiquée du paragraphe suivant.

### 1.1. Définition

Soit $A$ un anneau sans $p$-torsion et $p$-adiquement complet.

1.1.1. *Caractères et poids*. — Soit $\omega : \mathbf{Z}_p^\times \to \mu_{p-1} \subset \mathbf{Z}_p^\times$ la composition de la réduction modulo $p$ et du relèvement de Teichmüller. On a un isomorphisme $\mathbf{Z}_p^\times \simeq \mu_{p-1} \times (1 + p\mathbf{Z}_p)$ donné par $x \mapsto (\omega(x), x/\omega(x))$. Ainsi, le groupe $X$ des caractères continus de $\mathbf{Z}_p^\times$ à valeurs dans $\mathbf{Z}_p$ est-il naturellement isomorphe à $(\mathbf{Z}/(p-1)\mathbf{Z}) \times \mathbf{Z}_p$. On associe en effet à $(i, k)$ le caractère

$$\chi_{i,k} : \mathbf{Z}_p^\times \to \mathbf{Z}_p^\times, \quad x \mapsto \omega(x)^i (x/\omega(x))^k.$$

Lorsque $N = 1$ et $A = \mathbf{Q}_p$, Serre a donné dans [25, 1.4, b)] la définition suivante d'une forme modulaire $p$-adique.

*Définition 1.1.2.* — Soient $\chi \in X$ et $N \geqslant 1$ un entier premier à $p$. On dit qu'une série $f \in A[[q]]$ est une forme modulaire $p$-adique de poids $\chi$ pour le sous-groupe $\Gamma_1(N)$ s'il existe pour tout entier $n \geqslant 1$ une série $f_n \in A[[q]]$, $q$-développement d'une forme modulaire de poids $k_n$ pour $\Gamma_1(N)$, telle que

i) $f \equiv f_n \pmod{p^n}$;

ii) pour tout $x \in \mathbf{Z}_p^\times$, $\chi(x) \equiv x^{k_n} \pmod{p^n}$.

Toute forme modulaire usuelle de poids $k \in \mathbf{N}$ pour $\Gamma_1(N)$ est évidemment une forme modulaire $p$-adique de poids le caractère $x \mapsto x^k$. On dira aussi qu'une forme modulaire $p$-adique est *parabolique* si son terme constant est nul.

En utilisant la théorie des formes modulaires modulo $p$ de Swinnerton-Dyer (voir [26], ainsi que [24] et [25], 1.2), Serre prouve, dans le cas $N = 1$, que pour tout suite $(f_n)$ de formes modulaires de poids $k_n$ telle que $f_n \equiv f \pmod{p^n}$, la suite des caractères $x \mapsto x^{k_n}$ converge vers $\chi$ dans $X$ ([**25**], théorème 2, p. 202).

Il en déduit le résultat suivant (*loc. cit.*, corollaire 2, p. 204) :

PROPOSITION 1.1.3. — *Soit $(f_i = \sum a_n(f_i)q^n)$ une suite de formes modulaires p-adiques dans $\mathbf{Q}_p$, de niveau 1 et de poids $\chi_i \in X$. Supposons :*

i) *les coefficients $a_n(f_i) \in \mathbf{Q}_p$ pour $n \geqslant 1$ convergent uniformément vers $a_n \in \mathbf{Q}_p$ lorsque $i \to \infty$ ;*

ii) *les poids $\chi_i$ convergent vers un caractère $\chi \in X$ distinct du caractère unité $\mathbf{1}$.*

*Alors, le coefficient $a_0(f_i)$ a une limite $a_0 \in \mathbf{Q}_p$ lorsque $i \to \infty$ et la série*

$$f = \sum_{n=0}^{\infty} a_n q^n$$

*est une forme modulaire p-adique de poids $\chi$.*

Il en déduit aussi l'existence des « séries d'Eisenstein *p*-adiques », autrement dit l'existence d'une fonction continue, interpolation *p*-adique des valeurs de la fonction zêta de Riemann aux entiers négatifs.

*Exemple 1.1.4*. — Soit $\chi \in X$ un caractère *pair*, i.e. tel que $\chi(-1) = 1$. Pour toute suite $k_i$ d'entiers pairs qui converge vers $\chi$ dans $X$, et vers $+\infty$ dans $\mathbf{N}$, les séries

$$G_{k_i}^0 = \sum_{n=1}^{\infty} \sigma_{k_i-1}(n) q^n$$

(qui sont les séries d'Eisenstein $G_{k_i}$ privées de leur terme constant) convergent vers la série

$$G_\chi^{*,0} = \sum_{n=1}^{\infty} \bigg( \sum_{\substack{d|n \\ (d,p)=1}} \chi(d)/d \bigg) q^n.$$

Il résulte alors de la proposition 1.1.3 que les séries d'Eisenstein

$$G_{k_i} = \frac{1}{2}\zeta(1-k_i) + \sum_{n=1}^{\infty} \sigma_{k_i-1}(n) q^n$$

convergent vers une série

$$G_\chi = \frac{1}{2}\zeta(\mathbf{1}-\chi) + \sum_{n=1}^{\infty} \bigg( \sum_{\substack{d|n \\ (d,p)=1}} \chi(d)/d \bigg) q^n$$

que l'on appelle *série d'Eisenstein $p$-adique*. C'est une forme modulaire $p$-adique de poids $\chi$. Dans cette formule, $\mathbf{1} - \chi$ désigne le caractère $x \mapsto x/\chi(x)$ de $\mathbf{Z}_p^\times$, et $\zeta(1 - \chi)$ n'est autre que la valeur de la fonction zêta de Kubota–Leopoldt (la fonction zêta $p$-adique de $\mathbf{Q}$) en ce caractère.

1.1.5. *Définition alternative*. — Soit $K$ une extension finie de $\mathbf{Q}_p$ contenue dans $\mathbf{C}_p$ et $\mathfrak{O}$ son anneau d'entiers. Définissons sur l'espace vectoriel $\mathscr{M}_k(N; K)$ des formes modulaires de poids $k$ pour $\Gamma_1(N)$ à coefficients dans $A$ une norme $p$-adique en posant :

$$\|f\|_p = \sup_n |a_n(f)|_p, \quad f = \sum_{n=0}^\infty a_n(f) q^n.$$

On définit alors

$$\mathscr{M}(N; K) = \bigoplus_{i=0}^\infty \mathscr{M}_i(N; K) \subset K[[q]],$$

la somme étant directe, en vertu du principe de $q$-développement. On définit aussi

$$\mathscr{M}(N; \mathfrak{O}) = \mathscr{M}(N; K) \cap \mathfrak{O}[[q]] = \big\{ f \in \mathscr{M}(N; K) \mid \|f\|_p \leqslant 1 \big\}.$$

Alors, l'espace des formes modulaires $p$-adiques à coefficients dans $A = \mathfrak{O}$ ou $A = K$ n'est autre que le complété de $\mathscr{M}(N; A)$ pour la norme introduite.

1.1.6. Nous noterons $\overline{\mathscr{M}}(N; A)$ l'espace des formes modulaires $p$-adiques pour $\Gamma_1(N)$ à coefficients dans $A$, et $\overline{\mathscr{S}}(N; A)$ le sous-espaces des formes modulaires $p$-adiques paraboliques. Si $\chi \in X$, $\overline{\mathscr{M}}_\chi(N; A)$ et $\overline{\mathscr{S}}_\chi(N; A)$ désignerons les sous-modules des formes modulaires $p$-adiques de poids $\chi$. On notera $a_n(f)$ les coefficients d'une forme modulaire (éventuellement $p$-adique) $f$.

## 1.2. Opérateurs de Hecke

Pour simplifier, et parce que la théorie du § 2 est plus générale que celle que nous exposons ici, nous nous limitons désormais au cas $N = 1$.

1.2.1. *Opérateurs $T(\ell)$*. — Soit $f = \sum_{n=0}^\infty a_n(f) q^n$ une forme modulaire $p$-adique ; notons $\chi \in X$ son poids. Soit $(f_i)$ une suite de formes modulaires classiques de poids $k_i$ qui converge vers $f$, les $k_i$ convergeant vers $\chi$ dans $X$. Quitte à multiplier $f_i$ par une puissance assez grande de la série d'Eisenstein $E_{p-1}$, on peut supposer que $k_i$ tend vers $+\infty$ dans $\mathbf{R}$.

Soit $\ell$ un nombre premier. Alors, la suite

$$f_i | T(\ell) = \sum_{n=0}^\infty \Big( a_{\ell n}(f_i) + \ell^{k_i - 1} a_{n/\ell}(f_i) \Big) q^n$$

converge vers une forme modulaire $p$-adique quand $i \to \infty$, dont le développement de Fourier est

$$\sum_{n=0}^{\infty} \left(a_{\ell n}(f) + \chi(\ell)\ell^{-1}a_{n/\ell}(f)\right) q^n$$

si $\ell \neq p$, et est

$$\sum_{n=0}^{\infty} a_{pn}(f)q^n$$

si $\ell = p$. C'est par définition $f|T(\ell)$. Il est clair que $T(\ell)$ est un endomorphisme de $\overline{\mathscr{M}}(1;A)$. Quand $\ell \neq p$, c'est le prolongement par continuité de l'endomorphisme $T(\ell)$.

1.2.2. *L'opérateur $U$*. — C'est l'opérateur $T(p)$ ; il s'identifie en effet au prolongement par continuité de l'opérateur $U_p$ d'Atkin.

1.2.3. *L'opérateur $V$*. — Il est défini par la formule

$$f|V = \sum_{n=0}^{\infty} a_n(f)q^{np} = f(q^p).$$

Si $f$ est une forme modulaire classique, $f|V$ est une forme modulaire de niveau $Np$, si bien qu'il n'est pas automatique que $V$ préserve l'espace des formes modulaires $p$-adiques de niveau $N$. Cependant, soit $(f_i)$ une suite de formes modulaires classiques de poids $k_i$, avec $k_i \to +\infty$ dans $\mathbf{R}$, et $k_i \to \chi$ dans $X$, qui converge vers $f$. Alors,

$$f_i|V = p^{1-k_i}(f_i|T(p) - f_i|U)$$

est une forme modulaire $p$-adique de niveau 1 et poids $k_i$. Comme $f_i|V$ converge vers $f$, $f$ est elle-même une forme modulaire $p$-adique de poids $\chi$ et niveau 1.

1.2.4. *Exemple*. — Si $\chi \in X$ est un caractère pair, on a pour $\ell$ premier à $p$,

$$G_\chi^*|T(\ell) = (1 + \chi(\ell)/\ell)G_\chi^* \quad \text{et} \quad G_\chi^*|T(p) = G_\chi^*.$$

## 2. Point de vue modulaire

L'interprétation des formes modulaires classiques comme sections d'un certain faisceau inversible sur l'espace des modules des courbes elliptiques munies de structures de niveau adéquates s'est révélée extrêmement féconde dans la théorie des formes modulaires. Ce paragraphe est consacré à la description de la théorie géométrique des formes modulaires $p$-adiques, élaborée par Katz.

## 2.1. Trivialisations

Notons $N = N_0 p^r$, où $(p, N_0) = 1$.

*Définition 2.1.1.* — Soit $E$ une courbe elliptique sur un anneau $A$. Une *structure arithmétique pour* $\Gamma_1(N)$ *sur* $E$, appelée aussi $\Gamma_1(N)^{\mathrm{arith}}$-structure, est une inclusion
$$i : \mu_N \hookrightarrow E[N]$$
de schémas en groupes finis et plats sur $A$.

Si $r > 0$, c'est-à-dire si $p|N$, il importe de remarquer que l'existence d'une telle structure implique que les fibres géométriques de $E/A$ sont ordinaires. Si $N$ est inversible sur $A$, il s'agit d'une structure de niveau pour $\Gamma_1(N)$ usuelle, telle que définie par exemple dans [**17**].

2.1.2. Soit $N'$ un autre entier, multiple de $N$. Nous dirons que deux structures arithmétiques $i : \mu_N \hookrightarrow E[N]$ pour $\Gamma_1(N)$ et $i' : \mu_{N'} \hookrightarrow E[N']$ pour $\Gamma_1(N')$ sur une courbe elliptique $E/A$ sont *compatibles* si $i'|_{\mu_N} = i$.

*Définition 2.1.3.* — Soit $E$ une courbe elliptique sur un anneau $A$. Une *trivialisation* de $E$ est un isomorphisme de groupes formels
$$\varphi : \widehat{E} \xrightarrow{\sim} \widehat{\mathbf{G}_{\mathrm{m}}}.$$

Une condition nécessaire pour l'existence d'une trivialisation de $E$ est que les fibres géométriques de $E/A$ soient ordinaires. Réciproquement, une telle courbe elliptique $E/A$ possède des trivialisations après un changement de base fidèlement plat.

2.1.4. Nous dirons qu'une $\Gamma_1(N)^{\mathrm{arith}}$-structure $i$ sur $E/A$ est compatible avec une trivialisation $\varphi$ si l'application induite
$$\mu_{p^r} \hookrightarrow \widehat{E} \xrightarrow{\sim} \widehat{\mathbf{G}_{\mathrm{m}}}$$
est l'inclusion naturelle.

2.1.5. La donnée d'une trivialisation $\varphi$ est équivalente à la donnée, pour tout entier $\nu \geqslant 0$, d'une $\Gamma_1(p^\nu)^{\mathrm{arith}}$-structure, deux à deux compatibles, donnée qu'on pourrait nommer $\Gamma_1(p^\infty)^{\mathrm{arith}}$-structure.

## 2.2. Définition

Soit $A$ un anneau $p$-adique.

*Définition 2.2.1.* — Une fonction modulaire $p$-adique sur $A$ est la donnée pour tout quadruplet $(A', E/A', \varphi, i)$, où $A'$ est une $A$-algèbre $p$-adique, $E/A'$ une courbe elliptique, $\varphi$ une trivialisation de $E/A'$ et $i$ une structure arithmétique de niveau $N$ compatible à $\varphi$, d'un élément de $A$, compatibles à tout changement de base et ne dépendant que de la classe d'isomorphisme de $(E/A', \varphi, i)$.

On note $\mathbf{V}(N; A)$ la $A$-algèbre des fonctions modulaires $p$-adiques sur $A$.

2.2.2. *Construction géométrique*. — Le foncteur des algèbres $p$-adiques vers la catégorie des ensembles qui associe à tout anneau $p$-adique $A$ l'ensemble des classes d'isomorphismes de triplets $(E/A, \varphi, i)$ est représentable par un anneau $p$-adique qui n'est autre que $\mathbf{V}(N, \mathbf{Z}_p)$. La restriction de ce foncteur aux $A$-algèbres $p$-adiques est alors représenté par $\mathbf{V}(N; \mathbf{Z}_p) \hat{\otimes} A$.

Par définition, $\mathbf{V}(N; \mathbf{Z}_p) = \varprojlim_n \mathbf{V}(N; \mathbf{Z}/p^n\mathbf{Z})$, si bien qu'il suffit de décrire ce dernier anneau. Si $\nu \geqslant 0$, l'ensemble des classes d'isomorphisme de triplets $(A, E, i)$, où $A$ est une $\mathbf{Z}/p^n\mathbf{Z}$-algèbre, $E$ une courbe elliptique sur $A$ et $i$ une structure de niveau usuelle pour $\Gamma_1(Np^\nu)$ est représentable par la courbe modulaire affine $Y_1^*(Np^\nu)/(\mathbf{Z}/p^n\mathbf{Z})$. Soit $V_{n,m}$ l'anneau des fonctions de l'ouvert affine correspondant aux courbes elliptiques ordinaires, obtenu en enlevant l'ensemble fini des points supersinguliers. Alors, on a

$$\mathbf{V}(N; \mathbf{Z}_p) = \varprojlim_n \varinjlim_n \mathbf{V}_{n;m}.$$

2.2.3. On voit que cette construction ne dépend pas de la puissance exacte de $p$ qui divise $N$, autrement dit que $\mathbf{V}(N; \mathbf{Z}_p) = \mathbf{V}(N_0; \mathbf{Z}_p)$. C'est aussi clair sur la définition 2.2.1.

2.2.4. *Courbe de Tate*. — On connaît la courbe de Tate sur $\mathbf{Z}_p((q))$. Comme $\mathbf{Z}_p((q))$ n'étant pas $p$-adiquement complet, il convient de le compléter. On obtient l'anneau des séries de Laurent $\sum_{n \in \mathbf{Z}} a_n q^n$ telles que les sous-séries $\sum_{n \geqslant 0} a_n q^n$ et $\sum_{n \leqslant 0} a_n q^n$ convergent respectivement pour $|q|_p < 1$ et $|q|_p \geqslant 1$.

L'interprétation de $\mathrm{Tate}(q)$ comme quotient de $\mathbf{G}_\mathrm{m}$ par le sous-groupe engendré par $q$ fournit pour tout entier $N \geqslant 1$ une structure arithmétique pour $\Gamma_1(N)$ sur $\mathrm{Tate}(q)$, et elles sont compatibles. En particulier, on dispose d'un quadruplet canonique $(\widehat{\mathbf{Z}_p((q))}, \mathrm{Tate}(q), \varphi_\mathrm{can}, i_\mathrm{can})$ comme dans la définition 2.2.1 sur lequel on peut *évaluer* toute fonction modulaire $p$-adique.

Autrement, dit, on a défini pour tout fonction modulaire $p$-adique sur $\mathbf{Z}_p$ son développement de Fourier (ou $q$-développement), qui est un élément de $\widehat{\mathbf{Z}_p((q))}$. La définition s'étend à une fonction modulaire $p$-adique sur un anneau $p$-adique $A$.

*Définition 2.2.5.* — On dira qu'une fonction modulaire $p$-adique sur $A$ est holomorphe (resp. paraboliques) si son développement de Fourier appartient à $A[[q]]$ (resp. à $qA[[q]]$).

Nous noterons $\mathbf{W}(N; A)$ (resp. $\mathbf{W}_{\mathrm{cusp}}$) les idéaux de $\mathbf{V}(N; A)$ formés des fonctions modulaires $p$-adiques holomorphes (resp. paraboliques).

*Remarque 2.2.6.* — (a) Comme on peut trouver des formes différentielles sur les courbes modulaires sur $\mathbf{Z}/p^n\mathbf{Z}$ qui ne s'annulent pas hors des points singuliers, la construction fait en fait intervenir tous les poids.

(b) On pourrait aussi construire les fonctions modulaires holomorphes ou les paraboliques dans le style de 2.2.2. Pour plus de détails, voir [**6**], I.3.1, ainsi que les articles de Katz.

(c) Enfin, une mise en garde s'impose : Hida [**9**] et Gouvêa [**6**] ont des notations inversées l'un par rapport à l'autre. Nous avons suivi les notations de Hida ; Gouvêa note $\mathbf{W}$ ce que nous avons noté $\mathbf{V}$, et $\mathbf{V}$ ce que nous avons noté $\mathbf{W}$.

Le principe de $q$-développement reste valable pour les fonctions modulaires $p$-adiques :

THÉORÈME 2.2.7 (Principe de $q$-développement ([**13**], 1.4))

*Soit $A$ un anneau $p$-adique, l'application*

$$\mathbf{V}(N; A) \longrightarrow \widehat{A((q))}$$

*qui associe à une fonction modulaire $p$-adique sa $q$-développement est injective, de conoyau plat sur $A$.*

*En outre, si $A'$ est une $A$-algèbre $p$-adique contenant $A$, une fonction modulaire $p$-adique à coefficients dans $A'$ appartient en fait à $\mathbf{V}(N; A)$ si et seulement si sa $q$-développement est dans $\widehat{A((q))}$.*

La preuve de ce théorème repose sur le fait que l'application qui associe à une fonction modulaire $p$-adique à coefficients dans $A'$ son développement de Fourier est injective pour tout anneau $p$-adique $A'$, d'où la platitude. Grossièrement, la nullité du développement de Fourier implique que la fonction modulaire est nulle dans un voisinage de la courbe de Tate. Il s'agit d'en déduire la nullité de la fonction modulaire partout ; cela repose sur l'« irréductibilité » de l'espace des modules des courbes elliptiques trivialisées, voir [**14**] pour une démonstration dans cet esprit et [**13**] pour une preuve fondée sur la surjectivité de certaines représentations du groupe fondamental de la courbe modulaire en caractéristique $p$ privée de ses points supersinguliers.

### 2.3. Opérateurs de Hecke

2.3.1. *Opérateurs diamant*. — Soit $A$ un anneau $p$-adique. On va définir une action du groupe $Z = \mathbf{Z}_p^\times \times (\mathbf{Z}/N_0\mathbf{Z})^\times$ sur l'algèbre $\mathbf{V}(N_0; A)$ des fonctions modulaires $p$-adiques sur $A$. Soit en effet $(x, y) \in Z$. Si $f$ est une fonction modulaire $p$-adique, on définit une autre fonction modulaire $p$-adique $f|\langle x, y\rangle$ par la formule

$$(f|\langle x,y\rangle)(A', E/A', \varphi, i) = f(A', E/A', x^{-1}\varphi, yi).$$

2.3.2. *Opérateurs $T(\ell)$ pour $\ell \neq p$*. — Soit $\ell$ un nombre premier distinct de $p$. Soit $(A, E/A, \varphi, i)$ un quadruplet comme dans 2.2.1. Soit $H \subset E$ un sous-groupe d'ordre $\ell$ et notons $\pi : E \to E' = E/H$ le quotient de $E$ par $H$.

Comme $p \neq \ell$, $\pi$ induit un isomorphisme au niveau des groupes formels, si bien que la courbe elliptique $E'$ est munie d'une trivialisation naturelle $\varphi' = \varphi \circ \pi^{-1}$.

De plus, $\pi \circ i : \mu_N \to E'$ est une injection si $H$ n'est pas inclus dans le sous-groupe $i(\mu_N)$. On dispose dans ce cas d'une $\Gamma_1(N)^{\mathrm{arith}}$-structure sur $E'$. Si $\ell$ ne divise pas $N$, il existe (après un changement de base fidèlement plat) exactement $\ell + 1$ tels sous-groupes, tandis que si $\ell$ divise $N$, il n'y en a que $\ell$.

On définit alors l'opérateur $T(\ell)$ par

$$(f|T(\ell))(A, E, \varphi, i) = \ell^{-1} \sum_{\substack{H \hookrightarrow E \\ \#H = \ell \\ H \not\subset i(\mu_N)}} f(A, E/H, \varphi \circ \pi^{-1}, \pi \circ i).$$

(Par descente, c'est bien un élément de $A$.)

2.3.3. *Opérateur $U = T(p)$*. — Pour définir $T(p)$, le plus simple est de supposer, o.B.d.A., que $p|N$ et de remarquer, à l'aide du principe de $q$-développement, que la formule précédente définit un unique élément de $A$, si $A$ est sans $p$-torsion.

Il est aussi possible de définir un opérateur de Frobenius sur les fonctions modulaires $p$-adiques, en considérant le quotient d'une courbe elliptique trivialisée par son sous-groupe $\mu_p$ canonique. On peut alors interpréter l'opérateur $U$ comme la « trace » de Frob divisée par $p$. Enfin, mentionnons le fait que Frob : $\mathbf{V} \to \mathbf{V}$ est fini et plat de degré $p$, et que c'est un relèvement du morphisme de Frobenius sur $\mathbf{V}_1 = \mathbf{V}/p\mathbf{V}$ donné par l'élévation à la puissance $p$ ([**6**], II.2.9, p. 41).

2.3.4. *Action sur les $q$-développements*. — On peut prouver que les opérateurs que nous venons de définir vérifient les formules que l'on avait utilisées en 1.2 pour définir les opérateurs de Hecke sur les formes modulaires $p$-adiques du § 1.

En particulier, les idéaux des fonctions modulaires holomorphes (resp. paraboliques) sont préservées par les opérateurs de Hecke ainsi définis.

## 2.4. Comparaison avec la théorie du § 1, poids

*Définition 2.4.1*. — Soit $\chi : Z \to A^\times$ un caractère continu de $Z$. On dit qu'une fonction modulaire $p$-adique $f \in \mathbf{V}(N; A)$ est de poids $\chi$ si pour tout $(x, y) \in Z$, $f|\langle x, y\rangle = \chi(x, y)f$.

S'il existe $k \in \mathbf{Z}$ et un caractère $\varepsilon$ de $(\mathbf{Z}/N\mathbf{Z})^\times$ tels que $\chi(x, y) = x^k \varepsilon(y)$, on dira que $f$ est de poids $k$ et de Nebentypus $\varepsilon$.

On peut décomposer $\mathbf{V}(N_0; A)$ sous l'action de tout sous-groupe fini d'ordre premier à $p$ de $Z$. Cependant, la somme directe des sous-espaces correspondants à tous les caractères de $Z$ est une sous-algèbre de $\mathbf{V}$, *distincte de V* ! Voir [**6**], p. 18 pour plus de détails.

*Définition 2.4.2*. — Si $\chi \in X$ est un caractère de $\mathbf{Z}_p^\times$, on dira que $f \in \mathbf{V}(N; A)$ est de poids $\chi$ si pour tout $x \in \mathbf{Z}_p^\times$, $f|\langle x, 1\rangle = \chi(x)f$.

S'il existe $k \in \mathbf{Z}$ tel que $\chi(x) = x^k$, on dira que $f$ est de poids $k$.

Si $f$ est de poids $\chi = \chi_{(i,k)}$ pour $i \in \mathbf{Z}/(p-1)\mathbf{Z}$ et $k \in \mathbf{Z}$, on dira aussi que $f$ est de poids $(i, k)$.

Les formes modulaires $p$-adiques définies au § 1 sont exactement les fonctions modulaires de poids $\chi \in X$. Plus précisément :

THÉORÈME 2.4.3 (Katz, [**13**], Prop. A.1.6). — *Soient $A$ un anneau $p$-adiquement complet sans $p$-torsion et $\chi \in X$ un caractère continu de $\mathbf{Z}_p^\times$. Alors, une série $f \in A[[q]]$ est le $q$-développement d'une fonction modulaire $p$-adique (de niveau $N$) de poids $\chi$ si et seulement si c'est une forme modulaire $p$-adique de poids $\chi$ au sens de la définition 1.1.2*

Autrement dit, l'homomorphisme de $q$-développement induit un isomorphisme du $A$-module $\mathbf{V}_\chi(N; A)$ des fonctions modulaires $p$-adiques de poids $\chi$ sur $\overline{\mathscr{M}}_\chi(N; A)$. De plus, cet homomorphisme est compatible avec l'action des opérateurs de Hecke.

## 2.5. Anneau des congruences divisées

Soit $A$ un anneau $p$-adique de valuation discrète, et soit $K$ son corps des fractions.

*Définition 2.5.1*. — Le module des *congruences divisées* est la $A$-algèbre
$$\mathscr{D}(A) = A[[q]] \bigcap \bigoplus_k \mathscr{M}_k(N;K).$$

Cette algèbre est beaucoup plus grosse que l'algèbre des formes modulaires classiques sur $A$. Par exemple, $(E_{p-1} - 1)/p$ définit un élément de $\mathscr{D}(A)$.

PROPOSITION 2.5.2. — *On a une injection $\mathscr{D}(A) \subset \mathbf{V}(N;A)$.*

*Démonstration*. — En effet, soit $f \in \mathscr{D}(A)$. Par définition, il existe un entier $n \geqslant 0$ tel que $p^n f$ appartienne à l'algèbre des formes modulaires classiques. En particulier, $p^n f \in \mathbf{V}(N;A)$ et $f$ est un élément de $p$-torsion dans le quotient $\widehat{A((q))}/\mathbf{V}(N;A)$. D'après le théorème 2.2.7, ce quotient est plat sur $A$, soit sans $p$-torsion. Cela prouve que $f \in \mathbf{V}(N;A)$. □

## 3. Algèbre de Hecke ordinaire

### 3.1. Définitions de l'algèbre de Hecke $p$-adique

Les opérateurs de Hecke définis en 2.3 vérifient les mêmes relations que les opérateurs de Hecke agissant sur les formes modulaires classiques, de niveau divisant $p$. En particulier, ils commutent entre eux deux à deux.

Soient $A$ un anneau $p$-adique et $N$ un entier $\geqslant 1$. On sous-entendra souvent $(N;A)$ dans les notations, notant par exemple $\mathbf{W}$ pour $\mathbf{W}(N;A)$. On se restreint aux fonctions modulaires $p$-adiques holomorphes :

*Définition 3.1.1*. — L'algèbre de Hecke $\mathbf{T}(N;A)$ de $\mathbf{W} = \mathbf{W}(N;A)$ est la complétion pour la topologie compacte-ouverte[2] de la sous-algèbre commutative de $\operatorname{End}_{\mathbf{Z}_p} \mathbf{W}$ engendrée par les opérateurs $\mathbf{T}(\ell)$ pour $\ell \neq p$, $T(p) = U$ et les opérateurs diamant.

L'action des opérateurs diamant sur $\mathbf{W}$ définit une action continue du groupe $Z$ sur $\mathbf{T}$, et en particulier une action continue de $1 + p\mathbf{Z}_p$. Autrement dit, $\mathbf{T}$ est naturellement munie d'une structure de $\mathbf{Z}_p[[1+p\mathbf{Z}_p]]$-module. Rappelons que l'on note usuellement $\Lambda$ l'algèbre $\mathbf{Z}_p[[1+p\mathbf{Z}_p]]$ ; le choix d'un générateur de $1 + p\mathbf{Z}_p$, par exemple $\mathbf{u} = 1 + p$, fournit un isomorphisme $\Lambda \simeq \mathbf{Z}_p[[T]]$, par la formule $T \mapsto [\mathbf{u}] - 1$.

Les opérateurs de Hecke préservent l'espace $\mathbf{W}_{\mathrm{cusp}} = \mathbf{W}_{\mathrm{cusp}}(N;\mathbf{Z}_p)$ des formes paraboliques, d'où un quotient de $\mathbf{T}$ :

---

[2] Une base de cette topologie est fournie par les $\Omega_{K,\mathscr{U}}$, $K$ et $\mathscr{U}$ étant des parties de $\mathbf{W}$ respectivement compacte et ouverte, $\Omega_{K,\mathscr{U}}$ étant l'ensemble des endomorphismes $\varphi$ tels que $\varphi(K) \subset \mathscr{U}$.

*Définition 3.1.2.* — L'algèbre de Hecke $\mathbf{T}_0$ de $\mathbf{W}_{\text{cusp}}$ est le quotient de $\mathbf{T}$ par l'idéal anulateur de $\mathbf{W}_{\text{cusp}}$.

Enfin, un important théorème de Hida fournit un moyen de définir l'algèbre $\mathbf{T}_0$ qui n'utilise pas les formes modulaires $p$-adiques. Si $A$ est un anneau, soit $\mathfrak{h}_k(N; A)$ l'algèbre de Hecke usuelle agissant sur les formes paraboliques de poids $k$ et de niveau $N$. Si $p|N$, l'action des opérateurs de Hecke en niveau $Np$ sur une forme de niveau $N$ coïncide avec l'action des opérateurs de Hecke de niveau $N$. On a ainsi un homomorphisme naturel $\mathfrak{h}_k(Np; A) \to \mathfrak{h}_k(N; A)$.

*Définition 3.1.3.* — Si $A$ est un anneau $p$-adique, on définit
$$\mathfrak{h}_k(N_0 p^\infty; A) = \varprojlim_r \mathfrak{h}_k(N_0 p^r; A).$$

C'est de manière naturelle une $\Lambda$-algèbre et nous la munirons de la topologie limite projective des topologies $p$-adiques sur les $\mathfrak{h}_k(N_0 p^r; A)$.

On a vu 2.2.3 que les formes modulaires classiques (paraboliques) de niveau $N_0 p^r$ peuvent être considérées comme éléments de $\mathbf{W}_{\text{cusp}}(N; A)$ et que cette inclusion préserve les $q$-développements et est compatible aux opérateurs de Hecke (si $r \geqslant 1$). Cela fournit un homomorphisme continu et surjectif $\mathbf{T}_0 \to \mathfrak{h}_k(N_0 p^\infty; A)$.

Théorème 3.1.4 (Hida, [**10**], Théorème 4.1). — *Soit $A$ un anneau $p$-adique. La surjection naturelle $\mathbf{T}_0(N; A) \to \mathfrak{h}_k(Np^\infty; A)$ est un isomorphisme de $\Lambda$-algèbres.*

### 3.2. Dualités

3.2.1. Si $f$ est une forme modulaire classique et $T$ un élément de l'algèbre de Hecke, l'application $(f, T) \mapsto a_1(f|T)$ fournit un accouplement entre formes modulaires et éléments de l'algèbre de Hecke. Un des résultats importants de la théorie des formes modulaires classiques est que cet accouplement est parfait. Le but de ce paragraphe est d'étendre ce résultat aux formes modulaires $p$-adiques. Nous suivons la présentation de Gouvêa ([**6**], III.1), mais Hida a démontré des résultats similaires faisant intervenir une dualité de Pontryagin (*cf.* [**9**], §2).

On définit un accouplement
$$\mathbf{T}_0(N; A) \times \mathbf{W}_{\text{cusp}}(N; A) \longrightarrow A, \qquad (T, f) \mapsto a_1(f|T).$$

Si $f \in \mathbf{W}_{\text{cusp}}(N; A)$, l'application $T \mapsto a_1(f|T)$ définit un homomorphisme de $A$-modules $\Phi_f : \mathbf{T}_0(N; A) \to A$.

Le résultat fondamental est alors le suivant :

THÉORÈME 3.2.2. — *Soit $A$ l'anneau des entiers d'une extension finie $K$ de $\mathbf{Q}_p$. Soit $B$ une $A$-algèbre $p$-adique, munie de la topologie $p$-adique. L'application $f \mapsto \Phi_f$ induit un isomorphisme*

$$\mathbf{W}_{cusp}(N; B) \simeq \operatorname{Hom}_{A,cont}(\mathbf{T}_0(N; A), B).$$

(Rappelons que $\mathbf{T}_0 \subset \operatorname{End}_A \mathbf{W}_{\text{cusp}}(N; A)$ est munie de la topologie compacte-ouverte.)

*De plus, $\Phi_f$ est un homomorphisme de $A$-algèbres si et seulement si $f$ est une forme propre normalisée pour les opérateurs de Hecke et les opérateurs diamant.*

Il existe une variante du théorème précédent concernant les formes modulaires non paraboliques. Nous renvoyons le lecteur intéressé à l'article de Hida [9] et au livre de Gouvêa [6].

Enfin, comme la topologie $p$-adique sur $\mathbf{T}_0$ est strictement plus grossière que sa topologie limite projective, l'identité $\mathbf{T}_0 \to \mathbf{T}_0$ ne définit pas une forme modulaire $p$-adique à coefficients dans $\mathbf{T}_0$. Ce phénomène justifiera l'introduction des familles analytiques de formes modulaires $p$-adiques au § 4

### 3.3. L'algèbre de Hecke ordinaire

Commençons par un lemme général.

LEMME 3.3.1. — *Soient $A$ l'anneau des entiers d'une extension finie $K$ de $\mathbf{Q}_p$ et $H$ une $A$-algèbre finie. Pour tout $x \in H$, il existe un idempotent $e_x \in H$ vérifiant la propriété suivante : si $M$ est un $H$-module de type fini, $e_x M$ est le plus grand sous-module de $M$ sur lequel l'action de $x$ est un isomorphisme.*

*Démonstration.* — En fait, on va montrer que $e_x = \lim_{m \to \infty} x^{m!}$ convient. Il suffit en fait de démontrer le résultat pour la sous-$A$-algèbre $A[x] \subset \operatorname{End}_A(M)$. Cela permet de supposer que $H$ est commutative et finie sur $A$, en particulier complète pour la topologie $p$-adique.

Soit $\mathfrak{m}$ l'idéal maximal de $A$ et $k$ le corps résiduel $A/\mathfrak{m}$. Il suffit de prouver la convergence modulo $\mathfrak{m}$. En effet, si $n \geqslant 1$, $x^{m!}$ sera pour $m$ assez grand un idempotent de $H \otimes (A/\mathfrak{m}^n)$, d'où la convergence.

On décompose $x = s + n$ dans $H \otimes_A k$, où $s$ est un élément semi-simple, $n$ un élément nilpotent, et $xs = sx$. Or, pour $m$ assez grand $x^{m!} = s^{m!}$ est un projecteur puisque les valeurs propres de $s$ sont dans $\overline{k}^\times$, donc des racines de l'unité. □

Nous allons appliquer le lemme aux algèbres de Hecke (classique et $p$-adique) lorsque $x$ est l'opérateur $T(p)$.

*Définition 3.3.2.* — L'idempotent $e_k$ de $\mathfrak{h}_k(N; \mathbf{Z}_p)$ est l'idempotent attaché à $T(p)$ par le lemme 3.3.1.

On définit alors un idempotent $e$ de $\mathbf{T}_0 = \mathfrak{h}(Np^\infty; \mathbf{Z}_p)$ comme la limite de la suite des idempotents $e_k$ en niveaux $Np^r$.

Si $M$ est un $\mathbf{T}_0$-module, on appellera *partie ordinaire* de $M$ le sous-module $M^{\mathrm{ord}} = eM$. En particulier, *l'algèbre de Hecke ordinaire* est l'algèbre $\mathbf{T}_0^{\mathrm{ord}} = e\mathbf{T}_0$ facteur direct de $\mathbf{T}_0$.

Le lien avec les formes modulaires se fait de la façon suivante :

*Définition 3.3.3.* — On dit qu'une forme modulaire $f$ (classique ou $p$-adique) est ordinaire si $f|e = f$.

De manière concrète, comme $T(p)$ agit sur une forme propre par multiplication par $a_p$, on a la proposition :

PROPOSITION 3.3.4. — *Soit $A$ l'anneau des entiers d'une extension finie de $\mathbf{Q}_p$. Une forme modulaire $f$ $p$-adique (resp. classique), propre et normalisée (i.e. $a_1(f) = 1$) à coefficients dans $A$ est ordinaire si et seulement si $|a_p(f)|_p = 1$ (resp. et si de plus son niveau est divisible par $p$).*

Cette proposition explique le mot ordinaire : soit $E$ une courbe elliptique sur $\mathbf{Q}$ ayant bonne réduction en $p$. Alors, la courbe elliptique sur $\mathbf{F}_p$ réduction de $E$ est ordinaire si et seulement si le coefficient $a_p$ de sa fonction $L$ n'est pas multiple de $p$.

*Remarque 3.3.5.* — Soit $f$ une forme modulaire classique, propre et normalisée, de niveau $N_0$ premier à $p$, de poids $k \geqslant 2$ et de caractère $\chi$. Si $|a_p(f)|_p = 1$, on peut aisément calculer $f|e$. En effet, notons $U$ l'opérateur $T(p)$ de niveau $N_0p$, qui est distinct de l'opérateur $T(p)$ en niveau $N_0$ ; on a ainsi $U = T(p) + \chi(p)p^{k-1}V$. L'espace engendré par $f$ et $f|V$ est stable par $T(p)$, $U$ et $V$ et la matrice de l'endomorphisme $U$ dans la base $\{f, f|V\}$ est égale à

$$U = \begin{pmatrix} a_p & 1 \\ -\chi(p)p^{k-1} & 0 \end{pmatrix}.$$

Les valeurs propres $\alpha$ et $\beta$ de $U$ sont ainsi les solutions de l'équation

$$X^2 - a_p X + \chi(p)p^{k-1} = 0.$$

Les vecteurs propres de $U$ sont les
$$f_\alpha = f - \beta f|V \quad \text{et} \quad f_\beta = f - \alpha f|V.$$
On a $f = (\alpha f_\alpha - \beta f_\beta)/(\alpha - \beta)$, et quand $m \to \infty$,
$$U^{m!} f = \frac{\alpha^{1+m!}}{\alpha - \beta} f_\alpha - \frac{\beta^{1+m!}}{\alpha - \beta} f_\beta$$
converge vers $\alpha f_\alpha/(\alpha - \beta)$ si $\alpha$ est l'unique valeur propre qui est une unité $p$-adique. (Comme $k \geqslant 2$, $\alpha\beta \equiv 0 \pmod{p}$.) Autrement dit,
$$f|e = \frac{1}{2 - \alpha^{-1} a_p} \left( f(z) - \alpha^{-1} \chi(p) p^{k-1} f(pz) \right) = \frac{1}{2 - \alpha^{-1} a_p} f_0,$$
où $f_0$ est une forme modulaire ordinaire propre et normalisée. En particulier, on remarque que $f|e \neq f$, et donc que $f$ n'est pas ordinaire au sens de la définition 3.3.3.

Le théorème 3.2.2 s'étend en un résultat analogue de dualité entre formes modulaires paraboliques ordinaires et la partie ordinaire de l'algèbre de Hecke $\mathbf{T}_0$.

Nous pouvons maintenant énoncer le théorème fondamental de Hida concernant l'algèbre de Hecke ordinaire. Soient $r \geqslant 1$ et $k \geqslant 2$ deux entiers. On note le polynôme $\omega_{r,k} = (1+T)^{p^{r-1}} - \mathbf{u}^{p^{r-1}k}$ et posons $\Lambda_{r,k} = \Lambda/(\omega_{r,k})$.

THÉORÈME 3.3.6 (Hida). — *Soient $A$ un anneau $p$-adique et $N \geqslant 1$ un entier premier à $p$.*

i) *L'algèbre $\mathfrak{h}^{ord}(Np^\infty; A)$ est libre de rang fini sur $\Lambda$ ;*

ii) *On a un isomorphisme canonique*
$$\mathfrak{h}^{ord}(Np^\infty; A) \otimes_\Lambda \Lambda_{r,k} \simeq \mathfrak{h}^{ord}_k(Np^r).$$

Il existe deux preuves du point (i) de ce théorème. La première ([**8**, Théorème 3.1] et [**9**, Théorème 1.2]), valable pour $p \geqslant 5$, utilise les formes modulaires $p$-adiques et leur réduction modulo $p$. Elle utilise un résultat de N. Jochnowitz selon lequel l'image par l'opérateur $T(p)$ d'une forme modulaire modulo $p$ de poids $k + (p-1)$ pour $\Gamma_1$ est en fait de poids $k$. Il en résulte que la réduction modulo $p$ de l'espace des formes modulaires $p$-adiques ordinaires est un isomorphe à la somme directe des espaces de formes modulaires classiques modulo $p$ pour les poids $\leqslant p - 1$. Cette réduction est en particulier de dimension finie.

La seconde ([**28**], voir aussi [**11**], chapitre 7) fait usage de la théorie des formes modulaires $\Lambda$-adiques. Le point crucial est encore le fait que le rang de $\mathfrak{h}_k^{\mathrm{ord}}(N_0 p^r; \mathbf{Z}_p)$ est borné indépendamment du poids $k$ (*cf.* [**11**], chapitre 7, théorème 1, p. 202). Pour prouver ce fait, Wiles et Hida utilisent l'isomorphisme d'Eichler–Shimura qui relie formes modulaires et la cohomologie parabolique de $\Gamma_1(N)$ à valeurs dans certains $\mathbf{Z}[\Gamma_1(N)]$-modules. On peut alors établir que la partie ordinaire de cette cohomologie parabolique s'injecte dans un $\mathbf{F}_p$-espace vectoriel de dimension finie indépendant du poids.

## 4. Formes modulaires $\Lambda$-adiques

### 4.1. *Définition*

La notion de famille analytique de formes modulaires $p$-adiques a été introduite par Serre dans son article [**25**]. Soit $A$ l'anneau des entiers d'une extension finie de $\mathbf{Q}_p$; on note $\Lambda = A[[T]]$. On rappelle que $\mathbf{u} = 1 + p$.

*Définition 4.1.1.* — Soient $N$ un entier $\geqslant 1$ divisible par $p$ et $\psi$ un caractère de $(\mathbf{Z}/N\mathbf{Z})^\times$ On dit qu'une série $F(T;q) = \sum_{n=0}^\infty A_n(T) q^n \in \Lambda[[q]]$ est une forme modulaire $\Lambda$-adique de caractère $\psi$ s'il existe $r \geqslant 0$ tel pour presque tout $k \geqslant 2$, la série

$$F(\mathbf{u}^k - 1) = \sum_{n=0}^\infty A_n(\mathbf{u}^k - 1) q^n$$

est le développement de Fourier d'une forme modulaire de poids $k$, de niveau $N p^r$ et de caractère $\psi \omega^{-k}$.

Ainsi, se donner une forme modulaire $\Lambda$-adique revient à se donner une collection de formes modulaires classiques $f_k$, ainsi qu'une interpolation $p$-adique de leurs développements de Fourier par des éléments de $\Lambda$. Les formes modulaires classiques $F(\mathbf{u}^k - 1)$ sont appelées les *spécialisations* de $F$.

Nous dirons qu'une forme modulaire $\Lambda$-adique est parabolique (resp. ordinaire) si pour tout $k$ assez grand, sa spécialisation en poids $k$ est parabolique (resp. ordinaire). On notera $\mathsf{M}(N, \psi; \Lambda)$ (resp. $\mathsf{S}(N, \psi; \Lambda)$) les $\Lambda$-modules des formes modulaires $\Lambda$-adiques (resp. des formes modulaires $\Lambda$-adiques paraboliques) de niveau $N$ et caractère $\psi$.

4.1.2. *Séries d'Eisenstein.* — Soit $\psi_0$ un caractère pair de $(\mathbf{Z}/p\mathbf{Z})^\times$ à valeurs dans $A$. On démontre qu'il existe une série $E(\psi) \in \Lambda[[X]]$ si $\psi_0 \neq \mathbf{1}$ et $X^{-1}\Lambda[[X]]$ si $\psi_0 = \mathbf{1}$, telle que pour tout $r \geqslant 1$, tout caractère non trivial $\psi : (\mathbf{Z}/p^r\mathbf{Z})^\times \to \overline{\mathbf{Q}}^\times$ tel que $\psi|_{(\mathbf{Z}/p\mathbf{Z})^\times} = \psi_0$, et tout entier $k \geqslant 1$, on a

$$E(\psi_0)(\psi(\mathbf{u})\mathbf{u}^k - 1) = E_k(\psi \omega^{-k}) \in \mathscr{M}_k(\Gamma_1(p^r), \psi \omega^{-k}).$$

Ainsi, la série $E(\psi_0)$, *série d'Eisenstein $\Lambda$-adique,* interpole les séries d'Eisenstein $E_k(\psi\omega^{-k})$ usuelles. Le point crucial est l'interpolation du terme constant, c'est-à-dire des valeurs $L(1-k,\psi\omega^{-k})$, par une fonction analytique sur $\mathbf{Z}_p$ (Iwasawa). Pour plus de détails, *cf.* [**11**], chapitre VII, p. 198.

Une fois acquise l'existence de cette forme modulaire $\Lambda$-adique, on obtient de nombreuses formes modulaires $\Lambda$-adiques. Soit en effet $f$ une forme modulaire classique de poids $m$, de niveau $N$ divisible par $p$ et de caractère $\chi$ à valeurs dans $A$. Soit $\psi$ un caractère de $(\mathbf{Z}/N\mathbf{Z})^\times$ et $\psi_0$ sa restriction à $(\mathbf{Z}/p\mathbf{Z})^\times$ Le produit $fE(\psi_0)$ dans $\Lambda[[q]]$ n'est pas tout à fait une forme modulaire $\Lambda$-adique puisque $(fE(\psi_0))(\mathbf{u}^k - 1)$ est de poids $k + m$.

Posons alors $u_m = \psi(\mathbf{u})\mathbf{u}^{-m}$ et $v_m = u_m - 1$. Comme $|u_m|_p = 1$ et $|v_m|_p < 1$, la série
$$f * E(\psi_0) = (fE(\psi_0)(u_m X + v_m) = \sum_{n=0}^\infty a_n(u_m X + v_m)q^n$$
est bien définie ; c'est de plus une forme modulaire $\Lambda$-adique puisque si $k > m$,
$$\begin{aligned}f * E(\psi_0)(\mathbf{u}^k - 1) &= f(q)E(\psi_0)(\psi(\mathbf{u})\mathbf{u}^{-m}\mathbf{u}^k - 1) \\ &= f(q)E(\psi_0)(\psi(\mathbf{u})\,\mathbf{u}^{k-m} - 1) \\ &= f(q)E_{k-m}(\psi\omega^{m-k}).\end{aligned}$$
Ainsi, $f * E(\psi_0)$ est une forme $\Lambda$-adique de caractère $\psi\chi$. Si $f$ est parabolique, $f * E(\psi_0)$ est une forme $\Lambda$-adique parabolique.

4.1.3. *Opérateurs de Hecke.* — Les opérateurs diamant agissent naturellement sur les espaces des formes modulaires $\Lambda$-adiques, via leur caractère $\psi$. On peut faire agir les opérateurs de Hecke de manière relativement évidente : les formules usuelles donnant l'action des opérateurs de Hecke sur le $q$-développement d'une forme modulaire de poids $k$ (variable) s'interpolent facilement en une formule à coefficients dans $\Lambda$. Si $d$ est un entier premier à $p$, notons $s(d)$ l'unique élément de $\mathbf{Z}_p$ tel que $d = \omega(d)\mathbf{u}^{s(d)}$. Alors, la formule
$$\left(\sum_{m=0}^\infty A_m(T)q^m\right)|T(n) = \sum_{m=0}^\infty\left(\sum_{\substack{d|(m,n)\\(d,p)=1}}\chi(d)d^{-1}(1+T)^{s(d)}A_{mn/d^2}(T)\right)q^m$$
convient, puisque
$$\chi(d)d^{-1}(\mathbf{u}^k)^{s(d)} = \chi(d)d^{-1}(\mathbf{u}^{s(d)})^k = \chi(d)d^{-1}\omega(d)^{-k}.$$

On dira qu'une forme $\Lambda$-adique est propre si elle est propre pour l'action des opérateurs de Hecke que nous venons de définir, ce qui équivaut au fait que presque toutes ses spécialisations sont propres. (Un sens est clair, l'autre résulte du théorème de préparation de Weierstraß qui implique qu'un élément de $\Lambda$ est déterminé une fois qu'on le connait modulo une infinité d'idéaux premiers distincts.)

4.1.4. *Dualité*. — Comme nous n'avons énoncé le théorème de dualité que pour les formes paraboliques, nous nous restreignons désormais à ce cas. Une forme modulaire $\Lambda$-adique parabolique $F$ fournit par définition des formes modulaires $p$-adiques paraboliques $f_k$ par réduction modulo l'idéal $(1+T-\mathbf{u}^k)$ de $\Lambda$ (pour $k$ assez grand, disons $k \geqslant a$), d'où un homomorphisme de $A$-modules $\Phi_{f_k} : \mathbf{T}_0 \to A = \Lambda/(1+T-\mathbf{u}^k)$ continu pour la topologie $p$-adique. On en déduit un homomorphisme de $A$-modules

$$\Phi_F : \mathbf{T}_0 \longrightarrow \varprojlim_m \Lambda / \prod_{a \leqslant k \leqslant m} (1+T-\mathbf{u}^k) \simeq \Lambda.$$

Le dernier isomorphisme est dû au fait, conséquence du théorème de préparation de Weierstraß, qu'un élément de $\Lambda$ a un nombre fini de zéros dans le « disque unité ouvert » de $\mathbf{Z}_p$. Cet homomorphisme est continu lorsqu'on munit $\Lambda$ de sa topologie limite projective ; il n'est en revanche pas nécessairement continu lorsqu'on munit $\Lambda$ de la topologie $p$-adique. Ainsi, une forme modulaire $\Lambda$-adique *n'est pas* en général une forme modulaire à coefficients dans l'anneau $p$-adique $\Lambda$ !

Étudions maintenant l'action de $\Lambda$. Le caractère $\chi$ de $(\mathbf{Z}/N\mathbf{Z})^\times$ fournit un morphisme d'anneaux $\chi : \Lambda \to A$ et l'homomorphisme $\Phi_F : \mathbf{T}_0 \to \Lambda$ vérifie $\Phi_F((1+T)S) = \chi(\mathbf{u})\Phi_F(S)$ pour tout $S \in \mathbf{T}_0$.[3] Si l'on note $\Lambda_\chi$ le $\Lambda$-module $\Lambda \otimes_\chi A$, $\Phi_F$ est un homomorphisme de $\Lambda$-modules $\mathbf{T}_0 \to \Lambda_\chi$. On a ainsi le théorème suivant :

THÉORÈME 4.1.5. — *Le $\Lambda$-module des formes modulaires $\Lambda$-adiques paraboliques de niveau $N$ et de caractère $\chi$ est isomorphe à $\mathrm{Hom}_{\Lambda,cont}(\mathbf{T}_0, \Lambda_\chi)$, où $\mathbf{T}_0$ est muni de sa topologie ouverte-compacte, $\Lambda_\chi$ et $\Lambda_\chi$ de la topologie limite projective.*

On peut ainsi généraliser ce théorème en une définition des familles de formes modulaires à coefficients dans une $A$-algèbre $B$ topologique complète pour la topologie $p$-adique. Ce sont des $A$-homomorphismes $\Phi : \mathbf{T}_0 \to B$ tels que que pour tout homomorphisme d'algèbres $\nu : B \to A$, la composition $\nu \circ \Phi : \mathbf{T}_0 \to A$ est une forme modulaire $p$-adique à coefficients dans $A$, dite *spécialisation* de $\Phi$ via $\nu$.

---

[3] Attention, $T$ n'est pas un opérateur de Hecke dans cette formule...

4.1.6. Si $\varepsilon : \mathbf{Z}_p^\times \to A$ est un caractère d'ordre fini et $k$ un entier $\geqslant 2$, notons $\nu_{k,\varepsilon} : \Lambda \to A$ l'homomorphisme défini par $1+T \to \varepsilon(\mathbf{u})\mathbf{u}^k$. Il résulte alors du théorème précédent et de la description de l'algèbre de Hecke donnée dans le théorème 3.1.4 que $\nu \circ F$ est une forme modulaire $p$-adique de niveau $N$, de poids $k$, de caractère $\varepsilon\psi$ (de la définition résulte seulement le cas où $\varepsilon = \mathbf{1}$, avec en outre un nombre fini d'exceptions).

4.1.7. *Ordinarité*. — Soit $F$ une forme modulaire $\Lambda$-adique parabolique de niveau $N$ et de caractère $\psi$. Dire que $F$ est ordinaire revient à dire que presque toutes ses spécialisations $F_\nu$ sont ordinaires, c'est-à-dire dans $e\mathbf{W}_{\mathrm{cusp}}(N;A)$. Mais l'idempotent $e \in \mathbf{T}_0$ agit sur l'espace $\mathsf{M}(N,\psi;\Lambda)$. On constate ainsi que $(F|e)_\nu$, spécialisation par $\nu : (1+T) \mapsto \mathbf{u}^k$ de $F|e$ est égal à $F_\nu|e = F_\nu$ puisque $F_\nu$ est ordinaire. Le théorème de préparation de Weierstraß implique alors que $F|e = F$, c'est-à-dire que les deux notions d'ordinarité introduites (être dans l'image de $e$ et être ordinaire pour presque toute spécialisation) coïncident.

Enfin, si $F$ est propre et normalisée, on a vu dans la proposition 3.3.4 que $F_\nu$ est ordinaire si et seulement si $\nu(a_p(F)) = F(\mathbf{u}^k - 1)$ est une unité $p$-adique. Cela signifie que le terme constant de $a_p(F)$ est une unité $p$-adique, et donc que $a_p(F)$ est un élément inversible de $\Lambda$.

## 4.2. Représentations galoisiennes

4.2.1. Soit $f$ est une forme modulaire parabolique (classique), propre normalisée, pour $\Gamma_1(N)$, de caractère $\chi$ et de poids $k \geqslant 1$ à coefficients dans un corps de nombres $K$. Soit $\lambda$ une place de $K$, $p$ sa caractéristique résiduelle.[4] Alors, il existe une unique représentation

$$\rho : \mathrm{Gal}(\overline{\mathbf{Q}}/\mathbf{Q}) \to \mathrm{GL}(2, K_\lambda)$$

continue, absolument irréductible, non ramifiée hors de $Np$ et telle que pour tout nombre premier $\ell$ ne divisant pas $Np$,

$$\det(1 - \mathrm{Frob}_\ell X) = 1 - a_\ell(f)X + \chi(\ell)\ell^{k-1}X^2,$$

$\mathrm{Frob}_\ell$ étant un élément de Frobenius en la place $\ell$. L'existence d'une telle représentation (hormis la détermination précise de la ramification) est due à

---

[4] Je m'excuse de ce déséquilibre dans les notations mais préfère rester conforme aux us de la théorie d'Iwasawa.

Eichler–Shimura lorsque $k = 2$, à Deligne pour $k > 2$, et à Deligne–Serre lorsque $k = 1$.

4.2.2. Soit maintenant $F$ une forme modulaire $\Lambda$-adique de caractère $\psi$. Si $\varepsilon$ est un caractère d'ordre fini de $\mathbf{Z}_p^\times$, soit $\nu$ l'homomorphisme $\Lambda \to A[\chi]$ donné par $1 + T \to \chi(\mathbf{u})\mathbf{u}^k$ ; alors, $F_\nu = \nu \circ F$ est une forme modulaire de poids $k$ et de caractère $\psi\varepsilon$ (cf. 4.1.6). Si $F$ est propre normalisée, les $F_\nu$ aussi, d'où des représentations galoisiennes $\rho_\nu$ comme au paragraphe précédent. Il est raisonnable de se demander s'il existe une représentation galoisienne $\rho_F$ à valeurs dans $\mathrm{GL}(2,\Lambda)$ qui « interpole » les $\rho_\nu$ et Hida a prouvé dans [8] qu'il en est bien ainsi. Ce résultat a été reconsidéré par Wiles dans [28] (dans le cadre plus général des formes modulaires de Hilbert d'un corps totalement réel) et a fourni une nouvelle démonstration du théorème de Hida. En fait, il déduit l'existence de $\rho_F$ de l'existence d'une infinité des $\rho_\nu$. Au sein d'une famille $\Lambda$-adique de formes modulaires, cela permet de démontrer l'existence d'une représentation pour la spécialisation en poids $k \geqslant 2$ de l'existence des représentations galoisiennes en poids 2 pour une infinité de caractères. En particulier, comme toute forme modulaire ordinaire s'insère dans une famille analytique, on peut démontrer ainsi le résultat de Deligne, dans le cas ordinaire, en n'utilisant que le résultat d'Eichler–Shimura. On voit l'intérêt de la méthode pour le cas d'un corps totalement réel où l'on ne disposait pas du théorème de Deligne...

Théorème 4.2.3 (Hida). — *Soit $F$ une forme modulaire $\Lambda$-adique de niveau $N$ et de caractère $\psi$, propre normalisée et parabolique. Il existe alors une unique représentation*
$$\rho_F : \mathrm{Gal}(\overline{\mathbf{Q}}/\mathbf{Q}) \to \mathrm{GL}(2,\Lambda)$$
*vérifiant les propriétés suivantes :*

  i) *$\rho_F$ est continue pour la topologie p-adique, absolument irréductible ;*
  ii) *$\rho_F$ est non ramifiée hors de $Np$ ;*
  iii) *pour tout nombre premier $\ell$ ne divisant pas $Np$,*
$$\det(1 - \rho_F(\mathrm{Frob}_\ell)X) = 1 - a_\ell(F)X + \chi(\ell)(1+T)^{s(\ell)}\ell^{-1}X^2.$$

(On rappelle que $s(\ell)$ est l'unique élément de $\mathbf{Z}_p$ tel que $\ell = \omega(\ell)\mathbf{u}^{s(\ell)}$.)

### 4.3. *Indications sur la preuve de Wiles*

Nous allons expliquer la démonstration du résultat suivant, dû à Wiles. Compte-tenu de l'existence de représentations galoisiennes associées aux formes modulaires classiques, le théorème de Hida 4.2.3 en découle immédiatement.

PROPOSITION 4.3.1 (Wiles). — *Supposons que pour une infinité d'homomorphismes continus de $A$-algèbres $\nu : \Lambda \to \overline{\mathbf{Q}_p}$, notant $A_\nu$ la clôture intégrale de l'image de $\nu$, il existe une représentation $\rho_\nu : \mathrm{Gal}(\overline{\mathbf{Q}}/\mathbf{Q}) \to \mathrm{GL}(2, A_\nu)$ non ramifiée hors de $Np$ et telle que pour $\ell \nmid Np$,*

$$\det(1 - \rho_\nu(\mathrm{Frob}_\ell)X) = 1 - \nu(a_\ell(F))X + \chi(\ell)(1+T)^{s(\ell)}\ell^{-1}X^2.$$

*Soit $L$ le corps des fractions de $\Lambda$; il existe une représentation galoisienne $\rho : \mathrm{Gal}(\overline{\mathbf{Q}}/\mathbf{Q}) \to \mathrm{GL}(2, L)$ vérifiant les conditions du théorème 4.2.3.*

En fait, l'irréductibilité de la représentation $\rho_F$ se déduit du théorème analogue de Ribet (voir [**21**], théorème 2.3) pour les formes modulaires classiques. Le point crucial est donc la *construction* d'une représentation $\rho$ vérifiant la condition indiquée sur les Frobenius.

Typiquement, si $R$ est un anneau, et si $I$ et $J$ sont deux idéaux de $R$, disposant de représentations d'un groupe $G$ à valeurs dans $\mathrm{GL}(2, A/I)$ et $\mathrm{GL}(2, A/J)$ dont les réductions modulo $I + J$ sont isomorphes, on voudrait en déduire une représentation dans $\mathrm{GL}(2, A/(I \cap J))$. Formulée ainsi, la question est peut-être sans espoir. Ce qui est facile en revanche, c'est de déterminer le polynôme caractéristique des images de la représentation cherchée, et en particulier, les traces. Il faut néanmoins identifier les relations nécessairement vérifiées par les traces d'une représentation.

*Définition 4.3.2.* — Soient $G$ un groupe et $R$ un anneau topologiques. On appelle *pseudo-représentation* de $G$ à valeurs dans $R$ la donnée de trois fonctions continues $a : G \to R$, $d : G \to R$, $x : G \times G \to R$ et d'un élément $c \in G$ tel que $c^2 = 1$ vérifiant les relations suivantes :

  i) $a(gh) = a(g)a(h) + x(g,h)$ et $d(gh) = d(g)d(h) + x(h,g)$ ;
  ii) $x(gh, jk) = a(g)a(k)x(h,j) + a(k)d(h)x(g,j) + a(g)d(j)x(h,k)$
      $+ d(h)d(j)x(g,k)$ ;
  iii) $x(g,h)x(j,k) = x(g,k)x(j,h)$ ;
  iv) $a(1) = d(1) = d(c) = 1$ et $a(c) = -1$ ;
  v) $x(g,1) = x(1,g) = x(g,c) = x(c,g) = 0$.

La trace et le déterminant d'une pseudo-représentation sont les applications $g \mapsto a(g) + d(g)$ et $g \mapsto a(g)d(g) - x(g,g)$.

Si 2 est inversible dans $R$, on a les relations $a(g) = \frac{1}{2}(\mathrm{tr}(g) - \mathrm{tr}(gc))$, $d(g) = \frac{1}{2}(\mathrm{tr}(g) + \mathrm{tr}(gc))$ et $x(g,h) = a(gh) - a(g)a(h)$, si bien qu'une pseudo-représentation est déterminée par sa trace.

Soit $\rho : G \to \mathrm{GL}(2, R)$ une représentation continue et $c \in G$ un élément tel que $c^2 = 1$ et tel que $\det(\rho(c)) = -1$. Ainsi, $\rho(c)$ a pour valeurs propres $1$ et $-1$, et l'on peut supposer que $\rho(c) = \begin{pmatrix} -1 & 0 \\ 0 & 1 \end{pmatrix}$. Notons alors $\rho(g) = \begin{pmatrix} a(g) & b(g) \\ c(g) & d(g) \end{pmatrix}$. Les applications

$$a : g \mapsto a(g), \quad d : g \mapsto d(g), \quad \text{et} \quad x : (g, h) \mapsto a(gh) - a(g)a(h)$$

définissent une pseudo-représentation de $G$ à valeurs dans $R$.

4.3.3. Mettons nous maintenant sous les hypothèses de la proposition 4.3.1, les $\rho_\nu$ (pour $\nu$ dans un ensemble infini $N = \{\nu_1, \dots\}$ d'homomorphismes $\Lambda \to \overline{\mathbf{Q}_p}$) fournissent des pseudo-représentations $(a_\nu, d_\nu, x_\nu)$ du groupe de Galois $G$ de l'extension maximale non ramifiée hors de $Np$. En effet, il suffit de fixer un élément $c$ d'ordre 2 (conjugaison complexe). Le déterminant de $\rho_\nu(c)$ est alors $-1$ par la théorie classique des représentations galoisiennes associées aux formes modulaires. Notons $I_\nu \subset \Lambda$ le noyau de $\nu$, de sorte que la pseudo-représentation associée à $\rho_\nu$ est à valeurs dans $\Lambda/I_\nu$.

Leurs traces sont évidemment compatibles en $\mathrm{Frob}_\ell$, puisque ce sont pour tout $\nu$ la réduction de $a_\ell(F) \in \Lambda$ modulo $I_\nu$. En vertu du théorème de densité de Čebotarev, $\mathrm{tr}(\rho_\nu) = \mathrm{tr}(\rho_{\nu'})$ modulo $\nu + \nu'$, d'où une pseudo-représentation à valeurs dans $\Lambda/(I_\nu \cap I_{\nu'})$. On en déduit une pseudo-représentation de $G$ à valeurs dans $\varprojlim_r \Lambda/(I_{\nu_1} \cap \cdots \cap I_{\nu_r})$. En vertu du théorème de préparation de Weierstraß, cette dernière algèbre est isomorphe à $\Lambda$ (comme algèbre topologique), si bien que nous avons construit la pseudo-représentation $(a, d, x)$ attachée à la représentation galoisienne que nous cherchons.

4.3.4. Il reste à construire la représentation $\rho$. Il y a pour cela deux cas.

- si $x(g, h) = 0$ pour tout $h$, posons $\rho(g) = \begin{pmatrix} a(g) & \\ & d(g) \end{pmatrix}$. Alors, $g \mapsto \rho(g)$ convient.

- s'il existe $g_0$ et $h_0$ tels que $x(g_0, h_0) \neq 0$, posons $c(g) = x(g_0, g)$ et $b(g) = x(g, h_0)/x(g_0, h_0)$. Alors, $g \mapsto \rho(g) = \begin{pmatrix} a(g) & b(g) \\ c(g) & d(g) \end{pmatrix}$ convient.

4.3.5. Si les idéaux de hauteur 1 dans $\Lambda$ sont principaux, on peut même trouver $\rho$ à valeurs dans $\mathrm{GL}(2, \Lambda)$. Dans le deuxième cas ci-dessus, il suffit en effet de choisir un élément $\theta \in L^\times$ dont le diviseur est égal au diviseur de l'idéal engendré par les $x(g_0, h)$ quand $h$ varie. On pose alors $c(g) = x(g_0, g)/\theta$ et $b(g) = x(g, h_0)\theta/x(g_0, h_0)$.

### 4.4. *Spécialisation en poids* 1

4.4.1. Nous avons déjà remarqué que toutes les spécialisations en poids $k \geqslant 2$ d'une forme modulaire $\Lambda$-adique sont des ($q$-développements de) formes modulaires $p$-adiques. En poids 1, il se passe en revanche d'étranges phénomènes...

4.4.2. Mazur et Wiles ont étudié dans [20] la restriction de la représentation galoisienne de Hida à un groupe de décomposition en $p$, et en particulier la théorie de Hodge $p$-adique des représentations galoisiennes obtenues.

Soit $\Delta = q \prod_{n \geqslant 1}(1-q^n)^{24} = \sum_{n=1}^{\infty} \tau(n) q^n$ l'unique forme modulaire parabolique normalisée de poids 1 et de niveau 12, $\tau$ étant la fonction de Ramanujan. Fixons un nombre premier $p$ tel que $\tau(p)$ n'est pas multiple de $p$, de sorte que $\Delta$ est ordinaire en $p$. (D'après [7], excepté 2 411, tout nombre premier $p$ compris entre 11 et 65 063 convient.) Soit $\alpha$ l'unique racine de $X^2 - \tau(p)X + p^{11} = 0$ qui est une unité $p$-adique, de sorte que $\Delta|e = \Delta(z) - \alpha^{-11}\Delta(pz)$ est une forme modulaire ordinaire de niveau $p$, poids 12 et caractère trivial. Ainsi, il existe une unique forme modulaire $\Lambda$-adique de caractère $\omega^{12}$, $F_\Delta = \sum A_n(F_\Delta) q^n$, telle que $F_\Delta(\mathbf{u}^{12} - 1) = \Delta|e$. On dispose alors d'une représentation

$$\rho_\Delta : \operatorname{Gal}(\overline{\mathbf{Q}}/\mathbf{Q}) \to \operatorname{GL}(2, \Lambda) = \operatorname{GL}(2, \mathbf{Z}_p[[T]])$$

vérifiant les conditions du théorème 4.2.3.

PROPOSITION 4.4.3. — *Supposons que $p$ vérifie les conditions suivantes :*

i) *$p$ ne divise pas $\tau(p)$ ;*

ii) *$p \not\equiv 1 \pmod{11}$*[5] *;*

iii) *$p \notin \{11, 23, 691\}$.*

*Alors, l'image de $\rho_\Delta$ contient $\operatorname{SL}(2, \mathbf{Z}_p[[T]])$.*

La preuve repose sur le résultat de théorie des groupes suivant, dû à N. Boston ([20], Appendix, prop. 3) :

PROPOSITION 4.4.4 (Boston). — *Soient $\rho : \operatorname{Gal}(\overline{\mathbf{Q}}/\mathbf{Q}) \to \operatorname{GL}(2, \mathbf{Z}_p[[T]])$ une représentation continue, $\overline{\rho} : \operatorname{Gal}(\overline{\mathbf{Q}}/\mathbf{Q}) \to \operatorname{GL}(2, \mathbf{F}_p)$ sa réduction modulo l'idéal $(p, T)$. Soit $I_p \subset \operatorname{Gal}(\overline{\mathbf{Q}}/\mathbf{Q})$ un sous-groupe d'inertie en $p$. Supposons que $\rho$ satisfait les propriétés suivantes :*

i) *l'image $\rho(I_p)$ de l'inertie est incluse dans le sous-groupe des matrices de la forme $\begin{pmatrix} 1 & * \\ 0 & * \end{pmatrix}$ ;*

---

[5] Cette condition est malheureusement oubliée dans [20]...

ii) $\rho(I_p)$ contient une matrice de la forme $\begin{pmatrix} 1 & * \\ 0 & 1+\alpha T \end{pmatrix}$, pour un élément $\alpha \in \mathbf{Z}_p[[T]]^\times$ ;

iii) pour tout $b \in \mathbf{F}_p^\times$, $\overline{\rho}(I_p)$ contient une matrice de la forme $\begin{pmatrix} 1 & * \\ 0 & b \end{pmatrix}$ ;

iv) l'image de $\overline{\rho}$ contient $\mathrm{SL}(2, \mathbf{F}_p)$.

Alors, l'image de $\rho$ contient $\mathrm{SL}(2, \mathbf{Z}_p[[T]])$.

Nous renvoyons à l'article [20] pour les détails de la démonstration de cette proposition ; expliquons cependant pourquoi les hypothèses en sont vérifiées.

Comme $(p, T) = (p, 1+T-\mathbf{u}^{12})$, la représentation $\overline{\rho}_\Delta$ est la représentation modulo $p$ attachée à la forme modulaire $\Delta$. La condition (iv), i.e. le fait que l'image de $\overline{\rho}_\Delta$ contienne $\mathrm{SL}(2, \mathbf{F}_p)$ est un théorème de Swinnerton-Dyer ([26], corollaire, p. 31) ; c'est là qu'on utilise l'hypothèse $p \geqslant 11$ et $p \notin \{23, 691\}$.

D'autre part, la proposition 2, p. 247 de [20] implique la condition (i). En fait, il y est prouvé que la restriction de la représentation $\rho_\Delta$ à un sous-groupe de décomposition $G_p \simeq \mathrm{Gal}(\overline{\mathbf{Q}_p}/\mathbf{Q}_p)$ en $p$ est de la forme $\begin{pmatrix} \varepsilon_1 & * \\ 0 & \varepsilon_2 \end{pmatrix}$. De plus, ils prouvent que $\varepsilon_1$ est l'unique caractère non ramifié de $G_p$ tel que $\varepsilon_1(\mathrm{Frob}_p) = A_p(F_\Delta)$. Enfin, ils donnent au bas de la page 250 une formule pour $\varepsilon_2$ que nous allons expliciter. Si $\chi : G_p \to \mathbf{Z}_p^\times$ est le caractère cyclotomique donnant l'action de $G_p$ sur les racines de l'unité d'ordre une puissance de $p$, écrivons $\chi(g) = \omega(g)\chi_0(g)$, avec $\chi(g) \equiv 1 \pmod{p}$. On a alors

$$\varepsilon_2(g) = \varepsilon_1(g)^{-1}\omega(g)^{11}\chi_0(g)^{-1}[\chi_0(g)],$$

où $[\chi_0(g)] \in \Lambda \simeq \mathbf{Z}_p[[1+p\mathbf{Z}_p]]$.

Alors, il suffit de prendre pour $g$ un élément tel que $\chi_0(g) \equiv 1+p \pmod{p^2}$ pour obtenir une série $\varepsilon_2(g) \in \Lambda$ de la forme requise par la condition (ii).

Enfin, en vertu de la condition (i), la restriction à l'inertie de la représentation $\overline{\rho}_\Delta$ est de la forme $\begin{pmatrix} 1 & * \\ 0 & b \end{pmatrix}$. Le déterminant de $\overline{\rho}_\Delta$ est égale à $\omega^{11}$, son image contient donc $\mathbf{F}_p^\times$ puisque $p \not\equiv 1 \pmod{11}$, et le théorème de densité de Čebotarev implique que sa restriction à l'inertie en $p$ est elle-même surjective, d'où la condition (iii).

4.4.5. Considérons maintenant la représentation

$$\rho_{\Delta,1} : \mathrm{Gal}(\overline{\mathbf{Q}}/\mathbf{Q}) \to \mathrm{GL}(2, \mathbf{Z}_p)$$

spécialisation en poids 1 de $\rho_\Delta$. Il résulte du théorème de Mazur et Wiles que l'image de cette représentation galoisienne contient $\mathrm{SL}(2, \mathbf{Z}_p)$. Ceci entraîne que $\rho_{\Delta,1}$ *n'est pas* la représentation galoisienne associée à une forme modulaire classique de poids 1, puisque celles-ci sont d'image finie (cf. [4]).

En particulier, la spécialisation en poids 1 de la forme modulaire $\Lambda$-adique $F_\Delta$ n'est pas le développement de Fourier d'une forme modulaire classique de poids 1.

Dans certains cas, la situation est encore plus étrange :

PROPOSITION 4.4.6. — *Si $p \in \{13, 17, 19\}$, la restriction à un groupe de décomposition en $p$ de la représentation $p$-adique $\rho_{\Delta,1}$ n'est pas de Hodge–Tate.*

*Démonstration.* — En effet, si elle était de Hodge–Tate, la forme de la restriction de $\rho_\Delta$ à l'inertie $I_p$ montre qu'elle serait de poids $(0, 0)$. D'après un théorème de Sen [**22**], l'inertie agirait à travers un quotient fini. Or, la description de la restriction de $\rho_\Delta$ au sous-groupe de décomposition en $p$ montre que la restriction au sous-groupe d'inertie de $\mathrm{Gal}(\overline{\mathbf{Q}}_p/\mathbf{Q}_p(\zeta_p))$ est unipotente. Si l'inertie agit à travers un quotient fini, cela implique ainsi que la restriction de la représentation $\rho_\Delta$ à $\mathrm{Gal}(\overline{\mathbf{Q}}/\mathbf{Q}(\zeta_p))$ est non ramifiée. Autrement dit, $\rho_{\Delta,1}$ fournit une extension infinie non ramifiée de $\mathbf{Q}(\zeta_p)$. Comme un théorème d'Odlyzko borne le degré d'une telle extension si $p \leqslant 19$, nous avons établi une contradiction. $\square$

## 5. Formes modulaires surconvergentes

### 5.1. Motivation classique : la conjecture de Gouvêa–Mazur

Soit $N \geqslant 1$ un entier premier à $p$.

*Définition 5.1.1.* — La pente d'une forme modulaire propre pour $\Gamma_1(Np)$ à coefficients dans $\mathbf{C}_p$ est par définition la valuation $p$-adique de la valeur propre de $U_p$, normalisée par $v(p) = 1$.

Ainsi, les formes modulaires de pente nulle sont celles que nous avons appelées ordinaires.

Notons $\omega : (\mathbf{Z}/p\mathbf{Z})^\times \to \mathbf{Z}_p^\times$ le caractère de Teichmüller.

*Définition 5.1.2.* — Soient $\alpha \in \mathbf{Q}$, $k \geqslant 2$ et $\varepsilon$ un caractère de $(\mathbf{Z}/p\mathbf{Z})^\times$. Alors, $\mathbf{d}(k, \alpha, \varepsilon)$ est la dimension du $\mathbf{C}_p$-espace vectoriel engendré par les formes modulaires de poids $k$, de niveau $Np$, de caractère $\varepsilon\omega^{-k}$, propres et de pente $\alpha$.

Lorsque $\alpha = 0$, la structure de l'algèbre de Hecke $p$-adique ordinaire établie par Hida montre que $\mathbf{d}(k, 0, \varepsilon)$ est indépendant du poids $k$. Cela motive, si l'on veut, la conjecture suivante, due à Gouvêa et Mazur :

CONJECTURE 5.1.3 (Gouvêa–Mazur). — *Soit $\alpha \in \mathbf{Q}$ ; soient $k$ et $k'$ deux entiers $> \alpha + 1$. Alors, pour tout caractère $\varepsilon$, si $|k - k'|_p < p^{-\alpha}$,*

$$\mathbf{d}(k, \alpha, \varepsilon) = \mathbf{d}(k', \alpha, \varepsilon).$$

*En outre, si ces dimensions valent 1, les développements de Fourier des uniques formes modulaires normalisées de poids $k$ et $k'$, caractère $\varepsilon$ et pente $\alpha$ sont égaux modulo $p(k' - k)$.*

Wan et Buzzard ont établi une forme affaiblie de la première partie de cette conjecture, où $p^{-\alpha}$ est remplacé par $p^{-A\alpha^2 + B\alpha + C}$ pour certaines constantes ne dépendant que de $N$ et $p$ (avec $A > 0$). Sous l'hypothèse que cette congruence plus forte est vérifiée, la seconde partie est un théorème de Coleman [**3**].

Une idée naturelle est de considérer les espaces de formes modulaires de poids $k$ « en famille » et de déduire de propriétés de continuité $p$-adique de l'opérateur $U$ la continuité requise pour les dimensions ou les vecteurs propres.

Une telle approche est sans espoir sur l'espace des formes modulaires $p$-adiques introduit aux paragraphes précédents. En effet, si $A$ l'anneau des entiers d'une extension finie de $\mathbf{Q}_p$, tout élément $\lambda \in A$ tel que $|\lambda|_p < 1$ est valeur propre de l'opérateur $U$ sur l'espace des formes modulaires $p$-adiques de poids $k$ à coefficients dans $A$, l'espace propre correspondant étant de dimension *infinie* :

PROPOSITION 5.1.4. — *On a une décomposition*

$$\overline{\mathscr{M}}_\chi(N; A) = \operatorname{Im} \operatorname{Frob} \oplus \ker U,$$

*le noyau de $U$ est de dimension infinie et pour tout $\lambda$ dans l'idéal maximal de $A$, $\ker U$ et $\ker(U - \lambda)$ sont topologiquement isomorphes.*

*Idée de la preuve*. — On établit sur les développement de Fourier la relation

$$(\operatorname{Frob}(f)g)|U = f(g|U)$$

pour tout couple $(f, g)$ d'éléments de $\mathbf{W}(N; A)$. On en déduit une suite exacte, dans laquelle $\overline{\mathscr{M}}_k(N; A)$ est noté $\overline{\mathscr{M}}$ :

$$0 \to \overline{\mathscr{M}} \xrightarrow{\operatorname{Frob}} \overline{\mathscr{M}} \xrightarrow{1 - \operatorname{Frob} \circ U} \overline{\mathscr{M}} \xrightarrow{U} \overline{\mathscr{M}} \to 0.$$

Ainsi, $\ker U$ est égal au conoyau de Frob : $\overline{\mathscr{M}} \to \overline{\mathscr{M}}$, lequel est de dimension infinie.

D'autre part, si $f_0 \in \ker U$ et si $\lambda \in A$ vérifie $|\lambda|_p < 1$, la série

$$f_\lambda = f_0 + \lambda \operatorname{Frob}(f_0) + \lambda^2 \operatorname{Frob}^2(f_0) + \cdots$$

converge dans $\overline{\mathscr{M}}$ et vérifie

$$\begin{aligned} f_\lambda|U &= f_0|U + \lambda \operatorname{Frob}(f_0)|U + \lambda^2 \operatorname{Frob}^2(f_0)|U + \cdots \\ &= 0 + \lambda f_0 + \lambda^2 \operatorname{Frob}(f_0) + \cdots = \lambda f_\lambda. \end{aligned}$$

$\square$

Ainsi, pour avoir une bonne théorie spectrale, il faut imposer des conditions supplémentaires sur les formes modulaires $p$-adiques. Il faut en particulier éviter de pouvoir appliquer indéfiniment l'opérateur Frob à nos nouvelles formes modulaires. La condition convenable, introduite par Dwork dans [**5**], consiste à imposer que la forme modulaire $p$-adique, a priori évaluable aux courbes elliptiques ordinaires, se prolonge aux courbes elliptiques « pas trop » supersingulières.

### 5.2. Définition des formes modulaires surconvergentes

5.2.1. *Rappels sur l'invariant de Hasse.* — Sur une courbe elliptique $E/R$, $R$ étant une $\mathbf{F}_p$-algèbre, le morphisme de Frobenius absolu $F_{\mathrm{abs}} : E \to E$ induit un endomorphisme $p$-linéaire de $H^1(E, \mathscr{O}_E) = \omega_E^{-1}$. Si $\eta$ est une base de $H^1(E, \mathscr{O}_E)$, $F_{\mathrm{abs}}^*(\lambda\eta) = \lambda^p F^*(\eta)$, si bien que

$$F_{\mathrm{abs}}^*(\eta) \otimes \eta^{\otimes(-p)}$$

est un élément bien défini de $\omega_E^{\otimes(p-1)}$. C'est *l'invariant de Hasse* de $E/S$. Il définit une forme modulaire $H$ de poids $p-1$ et de niveau 1. Un calcul sur la courbe de Tate, dû à Deligne, montre que le développement de Fourier de $H$ est égal à 1.

Soit d'autre part $E_{p-1}$ la série d'Eisenstein de niveau 1 et de poids $p-1$ dont le développement de Fourier est

$$E_{p-1} = 1 - \frac{2(p-1)}{B_{p-1}} \sum_{n=1}^\infty \sigma_{p-2}(n) q^n.$$

Modulo $p$, ce $q$-développement se réduit à 1. Le fait que $H$ et $E_{p-1}$ aient mêmes développement de Fourier implique, en vertu du principe de $q$-développement, que $H \equiv E_{p-1} \pmod{p}$.

*Définition 5.2.2.* — Soient $A$ un anneau $p$-adique, $N \geqslant 1$ un entier premier à $p$ et $k$ un entier. Soit aussi $r \in A$.

Une forme modulaire $p$-adique de niveau $N$, de poids $k$ et de condition de croissance $r$ à coefficients dans $A$ est une règle qui associe à tout triplet $(E/B, i, Y)$ formé d'une courbe elliptique $E$ sur une $A$-algèbre $B$, d'une structure $\Gamma_1(N)^{\mathrm{arith}}$ $i : \mu_N \hookrightarrow E[N]$ sur $E/B$ et d'une section $Y \in \omega_E^{\otimes(1-p)}$ telle que $Y \cdot E_{p-1}(E/B) = r$, un élément $f(E/B, i, Y) \in \omega_E^{\otimes k}$ ne dépendant que de la classe d'isomorphisme du triplet $(E/B, i, Y)$ et de formation compatible au changement de base.

5.2.3. On peut évaluer une telle forme modulaire sur la courbe de Tate
$$(\text{Tate}(q)/A((q)), i_{\text{can}}, Y_{\text{can}}),$$
où
$$Y_{\text{can}} = r \cdot E_{p-1}(\text{Tate}(q))^{-1},$$
ce qui a un sens car la réduction modulo $p$ de la courbe de Tate étant ordinaire, $E_{p-1}(\text{Tate}(q))$ est une base de $\omega_{\text{Tate}(q)}^{p-1}$. On en déduit un développement de Fourier dans $A((q))$.

Nous noterons $\mathsf{M}_k(N, r; A)$ le $B$-modules des telles formes modulaires $p$-adiques dont le développement de Fourier appartient à $A[[q]]$ et $\mathsf{S}_k(N, r; A)$ le sous-$B$-modules de celles — dites paraboliques — dont le développement de Fourier appartient à $qA[[q]]$. Pour $r = 1$, on obtient les espaces $\overline{\mathscr{M}}_k(N; A)$ et $\overline{\mathscr{S}}_k(N; A)$ étudiés auparavant. Si $r$ n'est pas une unité, les formes modulaires $p$-adiques correspondantes seront dites *surconvergentes*.

Si $r_2 = rr_1$, on dispose d'une application naturelle
$$\mathsf{M}_k(N, r_2; A) \to \mathsf{M}_k(n, r_1; A)$$
en associant à une forme $r_2$ de condition de croissance $r_2$ la forme $f_1$ définie par $f_1(E/B, i, Y) = f_2(E/B, i, rY)$ dont la condition de croissance est $r_1$.

5.2.4. *Interprétation rigide-analytique*. — Supposons $N \geqslant 3$ dans ce paragraphe. Notons $A = W(\overline{\mathbf{F}_p})$ et considérons la courbe modulaire $X_1(N)_A$, que nous noterons $X$ pour simplifier, qui classifie les courbes elliptiques $E$ munies d'une injection de $\mu_N$ dans $E[N]$. Dans sa réduction modulo $p$, $X_{\overline{\mathbf{F}_p}}$, il y a un nombre fini, disons $s$, de points correspondant à des courbes elliptiques supersingulières : $\overline{x}_1, \ldots, \overline{x}_s$. Pour tout $i$, choisissons une coordonnée locale $t_i$ de $X$ au voisinage de $\overline{x}_i$, de sorte que l'anneau local complété de $X$ en $\overline{x}_i$ est isomorphe à $A[[t_i]]$. Cela identifie pour tout $i$ l'ensemble des points de $X(\mathbf{C}_p)$ dont la réduction modulo $p$ est $\overline{x}_i$ au disque ouvert de rayon 1 dans $\mathbf{C}_p$ défini par $|t|_p < 1$.

Si on regarde $X(\mathbf{C}_p)$ comme espace analytique $p$-adique (rigide), il est loisible de lui ôter ces « disques supersinguliers ». On obtient alors un espace analytique rigide $X_1$. Mais on peut enlever des disques plus petits, à savoir ceux définis par l'inégalité $|t_i|_p < r$, pour $r \in p^{\mathbf{Q}}$. On obtient alors un autre espace analytique rigide $X_r$ qui contient $X_1$.

De plus, sur $X(\mathbf{C}_p)$, on dispose d'un faisceau $\omega$, qui est en quelque sorte l'analytifié du faisceau correspondant sur $X$.

La comparaison des points de vues géométrie formelle et géométrie rigide montre qu'après tensorisation par $\mathbf{Q}_p$, l'espace $\mathsf{M}_k(N, r; A)$ des formes modulaires de poids $k$, de niveau $N$ et de condition de croissance $r$ s'identifie aux sections (pour la géométrie analytique rigide) de $\omega^{\otimes k}$ sur $X_r$.

### 5.3. Bases de Banach

Katz a déterminé dans [**12**] la structure de ces espaces de formes modulaires $p$-adiques avec condition de croissance. On a en effet le résultat :

PROPOSITION 5.3.1. — *Soit $N$ un entier $\geq 3$ et $k \geq 1$. Si $k = 1$, on suppose de plus $N \leq 11$. Soit $A$ un anneau $p$-adique et $r \in A$ un élément non diviseur de zéro. On a alors un isomorphisme*

$$\mathsf{M}_k(N, r; A) \simeq \varprojlim_n \left( \bigoplus_{j=0}^{\infty} \mathscr{M}_{k+j(p-1)}(N; A) \right) \otimes_A (A/p^n A) / (E_{p-1} - r).$$

Cette proposition est la clef d'une interprétation plus explicite des formes modulaires avec condition de croissance en termes de formes modulaires classiques. Dans toute la suite de ce paragraphe § 5, nous ferons l'hypothèse que $N \geq 3$ et $k \geq 1$, et que si $k = 1$, alors $N \leq 11$.

5.3.2. Katz prouve tout d'abord que pour tous $k$ et $\alpha \geq 0$, l'homomorphisme injectif induit par la multiplication par $E_{p-1}$

$$\cdot E_{p-1} : \mathscr{M}_{k+\alpha(p-1)}(N; \mathbf{Z}_p) \longrightarrow \mathscr{M}_{k+(\alpha+1)(p-1)}(N; \mathbf{Z}_p)$$

admet une section. Choisissons alors pour tout $k$ une telle section une fois pour toute, et notons $B_k(N, \alpha+1; \mathbf{Z}_p)$ son image. On a ainsi pour tout $\alpha \geq 0$ une somme directe

$$\mathscr{M}_{k+(\alpha+1)(p-1)}(N; \mathbf{Z}_p) \simeq E_{p-1} \mathscr{M}_{k+\alpha(p-1)}(N; \mathbf{Z}_p) \oplus B_k(N, \alpha+1; \mathbf{Z}_p).$$

Notons $B_k(N, \alpha+1; A) = B_k(N, \alpha+1; \mathbf{Z}_p) \otimes A$, on a alors un isomorphisme

$$\bigoplus_{\alpha=0}^{j} B_k(N, \alpha; A) \simeq \mathscr{M}_{k+j(p-1)}(N; A), \quad \sum_{\alpha=0}^{j} b_\alpha \mapsto \sum_{\alpha=0}^{j} E_{p-1}^{j-\alpha} b_\alpha.$$

5.3.3. Notons $B_k^{\mathrm{rig}}(N; A)$ le sous-$A$-module de $\prod_\alpha B_k(N, \alpha; A)$ formé des suites $(b_\alpha)$ telles que $b_\alpha$ tend vers 0 pour la topologie $p$-adique lorsque $\alpha \to \infty$, c'est-à-dire telles qu'il existe une suite d'entiers $(n_\alpha)$ tendant vers $+\infty$ tels que $b_\alpha \in p^{n_\alpha} B_k(N, \alpha; A)$. C'est un $A$-module complet pour la topologie $p$-adique.

PROPOSITION 5.3.4. — *L'inclusion naturelle de $B_k^{rig}(N;A)$ dans le complété $p$-adique de $\bigoplus_\alpha \mathscr{M}_k(N;A)$ induit un isomorphisme*

$$B_k^{rig}(N;A) \to \mathsf{M}_k(N,r;A), \quad \sum_{\alpha=0}^\infty b_\alpha \mapsto \sum_{\alpha=0}^\infty r^\alpha b_\alpha E_{p-1}^{-\alpha},$$

*cette dernière expression étant la forme modulaire avec condition de croissance $r$ dont la valeur sur le triplet $(E/A, i, Y)$ est donnée par la série convergente*

$$\sum_{\alpha=0}^\infty b_\alpha(E/A, i) Y^\alpha.$$

COROLLAIRE 5.3.5. — *Supposons que $r_2 = rr_1$ et que $r_2$ n'est pas diviseur de $0$ dans $A$. Alors, l'homomorphisme naturel $\mathsf{M}_k(N, r_2; A) \to \mathsf{M}_k(N, r_1; A)$ est injectif, et s'exprime dans $B_k^{rig}(N;A)$ par l'application*

$$\sum_{\alpha=0}^\infty b_\alpha \mapsto \sum_{\alpha=0}^\infty r^\alpha b_\alpha.$$

Ce résultat a plusieurs conséquences importantes. La première est le fait que le principe de $q$-développement reste valide pour les formes modulaires avec condition de croissance (sur un anneau sans $p$-torsion). De plus, avec la condition de croissance donnée par $r = 1$, une forme modulaire $p$-adique est divisible par une puissance de $p$ si et seulement si son développement de Fourier l'est.

D'autre part, il permet de donner un critère sur le développement d'une forme modulaire $p$-adique image de $\sum b_\alpha \in B_k^{rig}$ pour que cette forme soit surconvergente avec condition de croissance $r$ : il faut et il suffit que pour tout $\alpha$, $r^\alpha$ divise $b_\alpha$ dans $\mathscr{M}_{k+\alpha(p-1)}(N;A)$.

Il en découle aussi le résultat fondamental :

COROLLAIRE 5.3.6. — *Supposons que $A$ est un anneau $p$-adique de valuation discrète de corps des fractions $K$ de caractéristique $0$, et soient $r$, $r_1$ et $r_2 \in A$, $r$ n'étant pas une unité de $A$, tels que $r_2 = rr_1 \neq 0$. Alors, l'inclusion canonique*

$$\mathsf{M}_k(N, r_2; A) \otimes_A K \longrightarrow \mathsf{M}_k(N, r_1; A) \otimes_A K$$

*est un homomorphisme complètement continu d'espaces de Banach $p$-adiques.*

### 5.4. Opérateurs de Hecke

5.4.1. On dispose d'une action de $(\mathbf{Z}/N\mathbf{Z})^\times$ sur $\mathsf{M}_k(N,r;A)$ qui résulte de l'action de $(\mathbf{Z}/N\mathbf{Z})^\times$ sur la structure $\Gamma_1(N)^{\mathrm{arith}}$.

Pour les opérateurs $T(\ell)$, avec $\ell \neq p$, la théorie est essentiellement la même que celle que nous avons évoquée en 2.3.2.

La théorie du sous-groupe canonique d'une courbe elliptique pas trop supersingulière, due à Lubin ([**19**], voir aussi [**12**]) montre que l'endomorphisme Frob déjà étudié se prolonge en un homomorphisme

$$\text{Frob} : \mathsf{M}_k(N, r^p; A) \to r^{-k}\mathsf{M}_k(N, r; A) \cap \mathsf{M}_k(N, 1; A)$$

dès que $|r|_p > p^{-1/(p+1)}$.

PROPOSITION 5.4.2. — *Soit $K$ le corps des fractions de $A$. Si $|r| > p^{-1/(p+1)}$, l'homomorphisme*

$$\text{Frob} : \mathsf{M}_k(N, r^p; A) \otimes K \to \mathsf{M}_k(N, r; A) \otimes K$$

*est fini étale de rang $p$.*

Ainsi, on peut encore définir l'opérateur $U$ comme la trace de Frob définie par $p$, d'où un homomorphisme que l'on prouve vérifier

$$U : \mathsf{M}_k(N, r; A) \to \frac{1}{p}\mathsf{M}_k(N, r^p; A)$$

si $|r| > p^{-1/(p+1)}$.

5.4.3. *Théorie spectrale.* — Supposons que $r$ n'est pas inversible dans $A$ et $|r| > p^{-p/(p+1)}$. La proposition précédente et le corollaire 5.3.6 impliquent que l'opérateur $U$ induit un endomorphisme complètement continu de l'espace de Banach $p$-adique des formes modulaires $p$-adiques surconvergentes de niveau $N$, poids $k$ et condition de croissance $r$. La théorie spectrale de ces opérateurs est due à Serre [**23**] ; ils regroupent à la fois les opérateurs compacts (c'est presque la définition), d'où une théorie de Fredholm, et les opérateurs à trace (car les valeurs propres tendent vers 0 dans $K$ qui est ultramétrique).

En particulier, on dispose d'un déterminant de Fredholm $P(t) = \det(1 - tU)$ qui est un élément de $K[[t]]$ qui converge pour tout $t \in K$ et tel que $\lambda \neq 0$ est une valeur propre de $U$ si et seulement si $P(1/\lambda) = 0$, la dimension de l'espace caractéristique associé étant précisément la multiplicité de $\lambda^{-1}$ comme racine de $P$. Ainsi, après tensorisation par $K$, les formes modulaires $p$-adiques (avec condition de croissance $r \notin A^\times$) de pente $\alpha \in \mathbf{Q}$ forment un $K$-espace vectoriel de dimension finie. C'est cette dimension dont Wan prouve la continuité $p$-adique, laquelle, jointe au résultat suivant, implique la forme affaiblie de la conjecture de Gouvêa–Mazur.

THÉORÈME 5.4.4 (Coleman, [2]). — *Toute forme modulaire p-adique surconvergente de poids k dont la pente est $< k - 1$ est classique.*

Il convient de remarquer qu'une forme classique de poids $k$ a une pente $\leqslant k - 1$. Il suffit en effet de le faire pour une forme $f$ propre pour les opérateurs de Hecke et de caractère $\varepsilon$. Alors, si $\lambda_1$ et $\lambda_2$ sont les racines de l'équation

$$X^2 - a_p(f)X + \varepsilon(p)p^{k-1} = 0,$$

les deux formes

$$f_1 = f - \lambda_2 f|V \quad \text{et} \quad f_2 = f - \lambda_1 f|V$$

sont propres pour $U_p$ avec valeurs propres $\lambda_1$ et $\lambda_2$. La somme des pentes étant $k - 1$, chacune a une pente $\leqslant k - 1$.

5.4.5. *Retour sur l'ordinarité*. — La théorie précédente s'applique aux formes modulaires $p$-adiques ordinaire. Elle implique qu'elles sont surconvergentes, avec la condition de croissance $r$ tel que $|r| = p^{-p/(p+1)}$. De toutes façons, la théorie de Hida implique que les formes modulaires $p$-adiques ordinaires sont des formes modulaires classiques, ce qui est le cas de pente 0 du théorème précédent.

## Épilogue : le théorème de Buzzard–Taylor

Je ne peux terminer ce texte sans évoquer un résultat récent de Buzzard et Taylor [1] dans lequel ils construisent une forme modulaire de poids 1 à l'aide d'à peu près toutes les notions que nous avons abordées, et encore d'autres...

THÉORÈME A.1 (Buzzard–Taylor). — *Soient $p \geqslant 5$ et $K$ une extension finie de $\mathbf{Q}_p$ ; notons $A$ l'anneau des entiers de $K$ et $\mathfrak{p}$ l'idéal maximal de $A$. Soit $\rho : \mathrm{Gal}(\overline{\mathbf{Q}}/\mathbf{Q}) \to \mathrm{GL}(2, A)$ une représentation continue vérifiant les conditions suivantes :*

   i) *$\rho$ n'est ramifiée qu'en un nombre fini de places ;*

   ii) *la réduction de $\rho$ modulo $\mathfrak{p}$ est absolument irréductible et modulaire ;*

   iii) *$\rho$ est non ramifiée à $p$ et $\rho(\mathrm{Frob}_p)$ a des valeurs propres $\alpha$ et $\beta$ distinctes modulo $\mathfrak{p}$.*

*Alors, $\rho$ est modulaire : il existe une forme modulaire classique de poids 1, propre et parabolique, $f = \sum a_m(f) q^m$, un plongement du corps $\mathbf{Q}(a_m(f))$ dans $K$ tel que pour presque tout nombre premier $\ell$, $a_\ell(f) = \mathrm{tr}\,\rho(\mathrm{Frob}_\ell)$.*

A.2.  Ce résultat a deux conséquences immédiates. Tout d'abord, l'image d'une telle représentation $\rho$ est finie, ainsi que le prédit une conjecture de Fontaine et Mazur. En effet, Deligne et Serre [4] ont prouvé que les représentations galoisiennes associées aux formes de poids 1 sont uniques et d'image finie.

D'autre part, pour tout plongement de $K$ dans $\mathbf{C}$, la fonction $L$ d'Artin associée à $\rho$ est entière, puisque c'est le cas pour les fonctions $L$ des formes modulaires de poids 1.

A.3.  Tâchons de donner les grandes lignes de la démonstration.

Ils montrent tout d'abord l'existence de deux formes modulaires modulo $\mathfrak{p}$, $\overline{f}_\alpha$ et $\overline{f}_\beta$, de poids 2 et de niveau $Np$ ($N$ étant un certain entier premier à $p$) dont la représentation galoisienne associée est $\rho$ modulo $\mathfrak{p}$. De plus, $\overline{f}_\alpha$ et $\overline{f}_\beta$ sont propres, la valeur propre de $U$ étant respectivement $\alpha$ et $\beta$ (modulo $\mathfrak{p}$).

Les théorèmes de Wiles, Taylor–Wiles et Diamond impliquent alors qu'il existe deux formes $\Lambda$-adiques $F_\alpha$ et $F_\beta$ de niveau $N$ relevant respectivement $\overline{f}_\alpha$ et $\overline{f}_\beta$ dont les représentations galoisiennes associées se spécialisent sur $\rho$ en poids 1. (Remarquer qu'un relèvement de $\overline{f}_\alpha$ ou $\overline{f}_\beta$ en une forme de poids 2 sera ordinaire, donc s'étendra en une forme modulaire $\Lambda$-adique par la théorie de Hida.)

Le problème est que rien ne garantit que la spécialisation de $F_\alpha$ (ou de $F_\beta$) en poids $k = 1$ correspondra à une forme modulaire de poids 1, comme l'exemple de Mazur–Wiles le montre.

Les spécialisations $f_\alpha$ et $f_\beta$ en poids 1 sont tout de même des formes modulaires $p$-adiques. Mieux, comme elles sont ordinaires, elles sont surconvergentes, avec une condition de croissance $p^{p/(p+1)}$ (cf. 5.4.5).

Ils posent alors

$$f = \frac{\alpha f_\alpha - \beta f_\beta}{\alpha - \beta} \quad \text{et} \quad f' = \frac{f_\alpha - f_\beta}{\alpha - \beta}.$$

On a $f' = f|V$, si bien que $f'$ est surconvergente, de condition de croissance $p^{1/(p+1)}$. Ils arrivent alors à construire une section $g$ sur les disques supersinguliers $|t_i|_p \leqslant p^{-1/(1+p)}$ à l'aide de $f'$ telle que de plus, $f$ et $g$ coïncident sur les couronnes $p^{-p/(1+p)} \leqslant |t_i|_p \leqslant p^{-1/(1+p)}$, d'où une section rigide analytique de $\omega$ sur toute la courbe $X_1(p)_{\mathbf{C}_p}$. Le théorème dit « GAGA rigide » implique que cette section est en fait algébrique, d'où la forme modulaire cherchée.

ANTOINE CHAMBERT-LOIR, Institut de mathématiques de Jussieu, Boite 247, 4, place Jussieu, F-75252 Paris Cedex 05 • *E-mail* : `chambert@math.jussieu.fr`